\documentclass[journal,twoside]{IEEEtran}
\usepackage{cite}

\ifCLASSINFOpdf
 \usepackage[pdftex]{graphicx}
\else
\fi
\usepackage[cmex10]{amsmath}
\usepackage{array}

\usepackage{mdwmath}
\usepackage{mdwtab}
\usepackage[caption=false,font=footnotesize]{subfig}
\usepackage{stfloats}
\usepackage{url}

\usepackage{rotating}

\usepackage[american,fulldiode]{circuitikz} 
\usepackage{graphicx} 
\usepackage{multirow}
\usepackage{amsfonts}
\usepackage{bbm}

\usepackage{multirow}
\usepackage[section]{placeins}
\usepackage{epstopdf}

\usepackage{slashbox}
\usepackage{amssymb}

\usepackage{lipsum,amsmath,multicol}

\usepackage{multirow}

\usepackage{float}
\usepackage{tablefootnote}
\usepackage[multiple]{footmisc}

\makeatletter
\def\ps@headings{%
\def\@oddhead{\mbox{}\scriptsize\rightmark \hfil \thepage}%
\def\@evenhead{\scriptsize\thepage \hfil \leftmark\mbox{}}%
\def\@oddfoot{}%
\def\@evenfoot{}}
\def\BState{\State\hskip-\ALG@thistlm}
\makeatother

\makeatletter
\newcommand{\removelatexerror}{\let\@latex@error\@gobble}
\makeatother

\usepackage{enumitem}

\allowdisplaybreaks

\usepackage{varwidth}
\newcolumntype{M}[1]{>{\begin{varwidth}[t]{#1}}l<{\end{varwidth}}}

\usepackage{algorithmic}

\usepackage[linesnumbered,ruled,vlined]{algorithm2e}



\begin{document}
%
\title{A Novel Fully Informed Water Cycle Algorithm \\for Solving Optimal Power Flow Problems \\in Electric Grids}
%
%
%
\author{\IEEEauthorblockN{Alireza Barzegar,$^{\ast}$ } \IEEEauthorblockN{Ali Sadollah,$^{\dagger}$, and  \IEEEauthorblockN{Rong Su,$^{\ast}$}}%
\thanks{\hspace*{-15pt}$\ast$: School of Electrical and Electronic Engineering, Nanayang Technological University, Singapore, alireza001@e.ntu.edu.sg, rsu@ntu.edu.sg.}%
\thanks{\hspace*{-15pt}$\dagger$: University of Science and Culture, Tehran, Iran sadollah@usc.ac.ir.}}%

\maketitle

\begin{abstract}
Optimal power flow (OPF) is a key tool for planning and operations in energy grids. The line-flow constraints, generator loading effect, piece-wise cost functions, emission, and voltage quality cost make the optimization model non-convex and computationally cumbersome to solve. Metaheuristic techniques for solving the problem have emerged as a promising solution to solve the complex OPF problem.  Recently, the water cycle algorithm (WCA), a method inspired by the observation of the water cycle process and the surface run-off model was proposed for solving optimization problems. This paper proposes an improved version of WCA that uses the concept of sharing global and local information among individuals to improve the exploitation ability compared with the standard WCA. The so called fully informed WCA (FIWCA) is tested against standard WCA and other metaheuristic techniques studied in the literature on IEEE 30 and 57 bus systems for various scenarios. Comparison and discussion regarding the performance and reliability of the metaheuristics approaches studied in literature are discussed. The obtained optimization results show that the better performance of proposed FIWCA comparing with the WCA and other algorithms especially in term of stability performance over replications. Consequently, it emerges as a tool for solving OPF in a reliable and efficient manner.
\end{abstract}



%
\IEEEpeerreviewmaketitle

\vspace*{-6pt}
\section{Introduction}
\label{l:introduction}
%
%
%
%

\IEEEPARstart{T}{ypically}, an optimal power flow (OPF) problem determines the values of grid operating or planning decision variables by solving a non-convex and nonlinear optimization problem \cite{lesieutre2005convexity} with underlying physical and operating constraints. The non-convex and nonlinear terms in the OPF optimization model are due to power flow conditions, generator loading effects, voltage cost, piecewise cost function and emission cost. This makes the OPF problem computationally intensive to solve and, usually, determining the globally optimal solution is a time-consuming process.
In literature, the mathematical programming approaches such as linear programming \cite{al2008established}, quadratic programming \cite{burchett1984quadratically}, mixed integer programming \cite{barzegar2015intelligent}, nonlinear programming \cite{habibollahzadeh1989hydrothermal} and Newton-based approaches and interior point methods \cite{yan1999improving} have been developed and employed to solve the OPF problem. However, these approaches usually provide a feasible solution, rather than a global one \cite{maffei2018semantic}. To overcome the computational difficulty as well as to obtain a global solution, meta-heuristic approaches such as genetic algorithms (GA) \cite{goldberg1988genetic}, simulated annealing (SA) \cite{kirkpatrick1983optimization}, particle swarm optimization (PSO) \cite{eberhart1995new}, ant colony optimization (ACO) \cite{colorni1991distributed}, and others have been studied in the literature.

Owing to the ability of these meta-heuristic techniques in finding global or near-global optimum of complex constrained optimization problems, their applications in electrical networks and specially planning and scheduling problems such as power allocation, load management, economic dispatch, and the OPF have been attracted researchers in recent years. For instances, using the evolutionary optimization methods such as SA \cite{roa2003solution}, PSO \cite{abido2002optimal}, fuzzy GA \cite{gaing2004real}, biogeography based optimization (BBO) \cite{roy2010biogeography}, teaching-learning based optimization (TLBO) \cite{ghasemi2015improved}, efficient evolutionary algorithm (EEA) \cite{reddy2014faster}, glowworm swarm optimization (GSO) \cite{reddy2016optimal}, and hybrid differential evolution and harmony search algorithm (Hybrid DE-HA) \cite{reddy2018optimal}.

However, since these approaches use meta-heuristics, performance improvements in these algorithms are sought to make them more deterministic and optimal. More recently, the water cycle algorithm (WCA), an optimization algorithm inspired by the water cycle process and observation of how streams and rivers flow into the sea, has been implemented for finding better optimal solutions of power flow optimization problems \cite{eskandar2012water}. Typically the WCA simulates the surface run-off process, i.e., one of the main phases in the water cycle process observed in nature, as updating formulation for generating new individuals during iterative optimization process \cite{sadollah2015water}.

So far, many applications in different fields of research have utilized the efficiency of WCA for solving complex optimization problems in the literature. For instance, some modified versions of WCA have been implemented in different applications such as rough set theory \cite{jabbar2014water}, detecting optimum reactive power dispatch problems \cite{lenin2014water}, optimal operation of reservoir systems \cite{haddad2014application} and antenna array pattern synthesis \cite{guney2015quantized}. While improvements on WCA continue to evolve, its use in solving OPF problem has also not been explored.

First introduced by Eberhart and Kennedy \cite{eberhart1995new} in 1995, the swarm optimization methods play an important role as the nature-inspired meta-heuristic optimization tools. Each of these algorithms is mainly induced by a phenomenon in nature such as fish schooling or bird flocking. In the canonical PSO algorithm every individual particle learns to update its velocity using the influence of the particle with the best so far obtained solution as well as the best obtained solution obtained by itself.

\begin{figure}
\centering
\centering
\vspace{5mm}
\subfloat[Complete]{\includegraphics[totalheight=0.12\textheight]{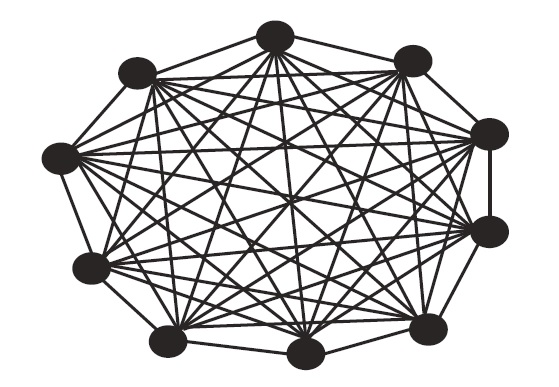}\label{topo_complete}}
\subfloat[Ring]{\includegraphics[totalheight=0.12\textheight]{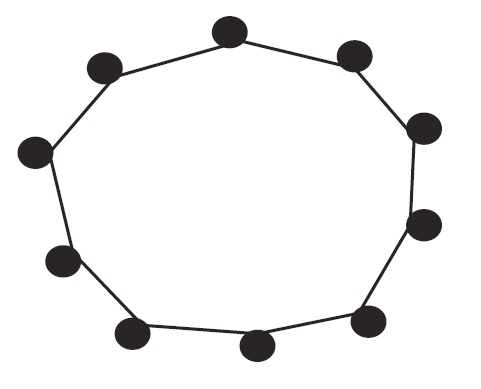}\label{topo_ring}}\\ \vspace{5mm}
\subfloat[Random]{\includegraphics[totalheight=0.12\textheight]{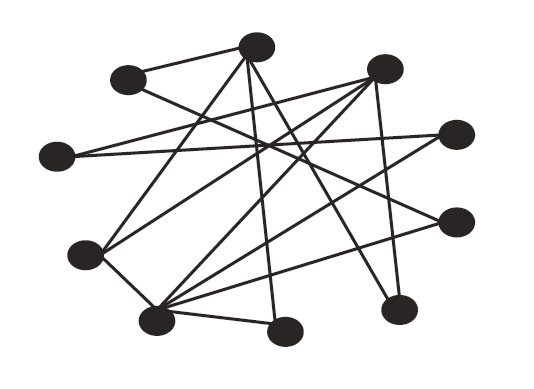}\label{topo_random}} \vspace{5mm}
\subfloat[Wheels]{\includegraphics[totalheight=0.12\textheight]{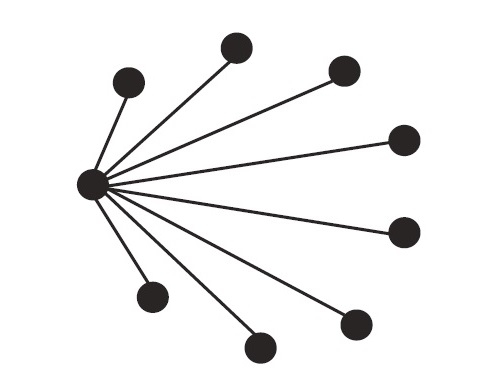}\label{topo_wheels}}\\ \vspace{5mm}
\subfloat[Star]{\includegraphics[totalheight=0.12\textheight]{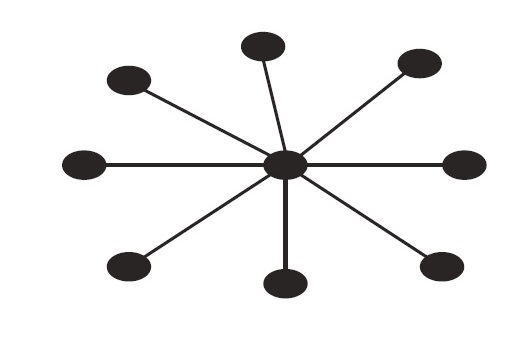}\label{topo_star}}
\subfloat[Von Neumann]{\includegraphics[totalheight=0.12\textheight]{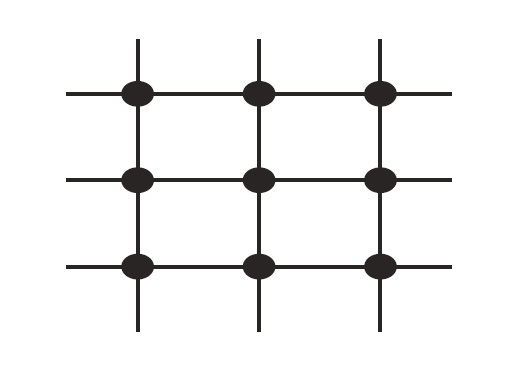}\label{topo_von_Neumann}}\\
\caption{Population Topologies in Particle Swarms \cite{lynn2018population}.}\vspace{5mm}
\label{Population_Topologies}
\end{figure}

The idea of population organization and structure in swarm algorithm was first presented in \cite{clerc2002particle}. Recent works show the significant impact of population structure and topology in improving the performance of these algorithms \cite{lynn2018population,payne2013complex,ali2016leveraged}. Figure \ref{Population_Topologies} demonstrates the most common population topologies in particle swarms and how each individual communicate with other members during the exploitation phase of swarm algorithms \cite{lynn2018population}.

Implementing the new communication topologies into the classical PSO results in new optimization new optimization meta-heuristics which aim to improve the performance of original algorithm by modifying the transferred information among individuals.
Mendes et al. proposed a fully informed particle swarm optimization (FIPSO) \cite{mendes2004fully} in which every particle's update strategy is influenced by all other members in the swarm. Some other algorithms also have been proposed based on the other communication topologies such as Firefly Algorithm (FA) \cite{yang2010nature,lukasik2009firefly} which benefits from ring/cluster topology or glowworm swarm optimization (GSO) \cite{krishnanand2009glowworm} in which every individual is influenced only by its better-performing neighbors within a specific radius. Therefore, considering the similarity between these algorithms and FIPSO, all these techniques can be represented as a family of Fully informed swarm (FIS) algorithms \cite{lukasik2014fully}.

Since communicating with all the swarm members may cause to considering redundant information, some FIS algorithms are proposed to utilize a more efficient of information from the neighbors instead of absorbing complete information of all the swarm members \cite{gao2015selectively,du2015adequate,du2017heterogeneous}.

Owing to the population-based properties of water cycle algorithm (WCA), communication topology among the individuals plays an important role in this optimization algorithm. The population individuals classified into three categories of sea, rivers and streams in which sea is a stream with the best solution and rivers are the next best solutions among the streams. Each stream may flows toward one of the rivers or sea, i.e. influenced by the information from only one of the rivers or sea, and the rivers may only flow towards the sea means that the each stream receives the information only from one of the rivers and/or sea and the rivers only communicate with the sea during the exploitation phase.

This investigation proposes an improved WCA that uses the idea of sharing global information among individuals during the exploitation phase, i.e. every stream benefits from the information not only from its allocated river/sea but also from all other rivers/sea during the exploration phase. Consequently, the proposed approach provides diversity in the selection of solution and eliminates the possibility of getting trapped at local optima, thereby enhancing the efficiency and accuracy of the algorithm significantly.

In the original WCA populations are clustered into a few sets, so we utilize this feature in the fully-informed algorithm and instead of using a complete communication topology among all the streams, only consider the information from the leader of each cluster, i.e. rivers or sea which leads to a more efficient algorithm in term of computation complexity while obtaining a higher performance toward finding the global solution.

The main contributions of this investigation are: (i) proposing a novel WCA with the concept of exchanging global and local information among individuals in the population, i.e. fully informed water cycle algorithm (FIWCA), (ii) solving the OPF problem considering nonlinearities introduced due to generator loading effect, line-flows, piece-wise cost, voltage quality cost and emission cost, and (iii) comparing the performance of the FIWCA with other meta-heuristic techniques including the traditional WCA for IEEE 30-bus, 57-bus test case, and renewable resources integrated 30-bus test case. Looking at obtained optimization results, the proposed FIWCA emerges as a promising technique for solving OPF with nonlinear cost/constraints and has the ability to reach the global optimum within a more reasonable time than existing techniques.

This paper is organized into 5 sections as given follows. The optimization model of the OPF and the different components are explained in Section 2. The proposed fully informed water cycle algorithm (FIWCA) is presented in Section 3 along with standard WCA. The optimization results, comparisons and discussions are provided in Section 4 for two reported test systems under different scenarios. Finally, conclusions and future research are given in Section 5.

\section{Problem Formulation}\label{sec2}

Consider an n-bus electric power system in which $\mathcal{N}=\{1,2,\ldots,n\}$ shows the set of all buses and $\mathcal{L}$ denotes the set of all lines where $(l,m)\in\mathcal{L}$ shows the line connecting bus $l$ to bus $m$. The optimal power flow (OPF) problem tries to minimize an objective function, $J(x,u)$, (e.g., generation cost) subject to several equality and inequality constraints, where $x$ is the state vector of dependent variables including the slack bus real power output, $P_{G1}$, load bus voltages, $V_{Dk}$, generators reactive power outputs, $Q_{Gk}$, and the transmission lines apparent power, $S_{lm}$. Vector $u$ shows the set of independent variables that contains the generators voltage $V_{Gk}$, generators real power outputs, $P_{Gk}$, except for the slack bus, $P_{G1}$, transformer tap settings, $T_{lm}$, and shunt VAR compensations, $Q_{Ck}$. The OPF problem deals with two groups of constraints including equality and inequality constraints given in the following subsections.

\subsection{Equality Constraints}
The equality constraints representing typical load flow equations are given as follows:
\begin{subequations}
\label{eq:PF_eq_BI_polar}
\begin{align}
\label{eq:PF_eq_BI_polar_P}
\nonumber
&\ P_{Gk} - P_{Dk} = |V_{k}|\sum_{i=1}^{n} |V_{i}|\Big(\mathbf{G}_{ik} \big(\cos({\theta}_i - {\theta}_k) \big)\\
& \qquad\qquad\qquad + \mathbf{B}_{ik} \big( \sin({\theta}_i - {\theta}_k) \big) \Big) & \forall k\in\mathcal{N}, \\
\nonumber
\label{eq:PF_eq_BI_polar_Q}
&\ Q_{Gk} - Q_{Dk} = |V_{k}|\sum_{i=1}^{n} |V_{i}|\Big(\mathbf{G}_{ik} \big(\sin({\theta}_i - {\theta}_k) \big)\\
& \qquad\qquad\qquad - \mathbf{B}_{ik} \big( \cos({\theta}_i - {\theta}_k) \big) \Big)  & \forall k\in\mathcal{N}.
\end{align}
\end{subequations}
where $P_{Gk}$ is the active power generation and $Q_{Gk}$ is the reactive power generation at bus \emph{k}-th and NB is the number of buses. $G_{ik}$, $B_{ij}$, and  are conductance, susceptance, and phase difference of voltages between bus \emph{k}-th and bus \emph{i}-th, respectively. The active and reactive load demands at bus \emph{k}-th are represented by $P_{Dk}$ and $Q_{Dk}$, respectively.

\subsection{Inequality constraints}
Inequality constraints include a set of system operating constraints, which are given as follows.  Note that the superscripts ``$min$'' and ``$max$'' denotes the lower and upper limits of the intended variable.
	
\subsubsection{Generator constraints}
the generator real and reactive power outputs and generator voltage are bounded as given follows:

\begin{equation}
\begin{aligned}
& P_{Gk}^{min} \leq P_{Gk} \leq P_{Gk}^{max}&     \forall k\in\mathcal{N},
\end{aligned}
\end{equation}
\begin{equation}
\begin{aligned}
& Q_{Gk}^{min} \leq Q_{Gk} \leq Q_{Gk}^{max} & \forall k\in\mathcal{N},
\end{aligned}
\end{equation}
\begin{equation}
\begin{aligned}
& \big|V_{Gk}^{min}\big| \leq |V_{Gk}| \leq \big|V_{Gk}^{max}\big| &  \forall k\in\mathcal{N}.
\end{aligned}
\end{equation}

\subsubsection{Transformer constraints}

The transformer tap settings are restricted by their lower and upper limits as given follows

\begin{equation}
\begin{aligned}
& T_{lm}^{min} \leq T_{lm} \leq T_{lm}^{max}  &  \forall (l,m)\in\mathcal{L}.
\end{aligned}
\end{equation}

\subsubsection{Shunt VAR compensator constraints}

The shunt VAR are bounded by their limits defined as provided follows:
\begin{equation}
\begin{aligned}
& Q_{Ck}^{min} \leq Q_{Ck} \leq Q_{Ck}^{max} & \forall k\in\mathcal{N}.
\end{aligned}
\end{equation}

\subsubsection{Security constraints}

The security limits include a set of constraints on load buses voltages and transmission lines powers given as follows
\begin{equation}
\begin{aligned}
& \big|V_{Lk}^{min}\big| \leq |V_{Lk}| \leq \big|V_{Lk}^{max}\big| &  \forall k\in\mathcal{N},
\end{aligned}
\end{equation}
\begin{equation}
\begin{aligned}
&\left| S_{lm} \right| \leq S_{lm}^{max} & \forall (l,m)\in\mathcal{L}.
\end{aligned}
\end{equation}
\vspace{1em}

The intended OPF problem is then
\begin{subequations}
\label{eq:opf_FIWCA}
\begin{align}
& \min \quad {J} \\
&\text{subject to}\\
\label{eq:opf_ch2_P}
& P_{Gk}^{min} \leq P_{Gk} \leq P_{Gk}^{max}&     \forall k\in\mathcal{N},\\
\label{eq:opf_ch2_Q}
& Q_{Gk}^{min} \leq Q_{Gk} \leq Q_{Gk}^{max} & \forall k\in\mathcal{N},\\
\label{eq:opf_ch2_V2}
& \big|V_k^{min}\big| \leq |V_{k}| \leq \big|V_k^{max}\big| &  \forall k\in\mathcal{N},\\
\label{eq:opf_ch2_S}
&\left| S_{lm} \right| \leq S_{lm}^{max} & \forall (l,m)\in\mathcal{L},\\
\label{eq:opf_ch2_Slm_Plm_Qlm}
&\left( S_{lm} \right)^2 = \left( P_{lm} \right)^2 + \left( Q_{lm} \right)^2 & \forall (l,m)\in\mathcal{L},\\
\label{eq:opf_ch2_PF_Eq_P}
\nonumber
&\ P_{Gk} - P_{Dk} = |V_{k}|\sum_{i=1}^{n} |V_{i}|\Big(\mathbf{G}_{ik} \big(\cos({\theta}_i - {\theta}_k) \big)\\
& \qquad\qquad\qquad + \mathbf{B}_{ik} \big( \sin({\theta}_i - {\theta}_k) \big) \Big) & \forall k\in\mathcal{N},\\
\label{eq:opf_ch2_PF_Eq_Q}
\nonumber
& Q_{Gk} - Q_{Dk} = |V_{k}|\sum_{i=1}^{n} |V_{i}|\Big(\mathbf{G}_{ik} \big(\sin({\theta}_i - {\theta}_k) \big)\\
& \qquad\qquad\qquad - \mathbf{B}_{ik} \big( \cos({\theta}_i - {\theta}_k) \big) \Big)  & \forall k\in\mathcal{N}.\\
\label{eq:opf_FIWCA_QC}
& Q_{Ck}^{min} \leq Q_{Ck} \leq Q_{Ck}^{max} & \forall k\in\mathcal{N},\\
\label{eq:opf_FIWCA_T}
& T_{lm}^{min} \leq T_{lm} \leq T_{lm}^{max}  &  \forall (l,m)\in\mathcal{L}.
\end{align}
\end{subequations}

The set of independent variables, i.e. control variables $u$, contains the generator voltage $V_{Gk}$, generator active power outputs $P_{Gk}$, except for the slack bus $P_{G1}$, transformer tap settings $T_{lm}$, and shunt VAR compensations $Q_{Ck}$.

\section{Water Cycle Algorithm Using the Concept of Global Information}

In this section, the standard Water Cycle Algorithm (WCA) and an improved WCA are described in details. Subsequently, Fully informed WCA, denoted as FIWCA, is proposed by equipping the standard WCA with the concepts of moving strategies toward the best solutions. In order to explain the FIWCA, first, the whole processes of WCA are given in the following section.

\subsection{Standard Water Cycle Algorithm}

The water cycle algorithm (WCA), as a meta-heuristic optimization technique, which emulates the natural water flow process first has been presented in \cite{eskandar2012water}. The algorithm initiates with the rain or precipitation phenomena by generation of a random population of design variables, or the streams population between lower and upper bounds. Then, the best stream, i.e. the stream with the minimum objective function value (for minimization problems) is selected as the sea.

Afterward, a set of streams with the closest objective values to the best objective value are selected as the rivers. Note that it is assumed that the remaining streams move toward the rivers and sea. Now let's see the mathematical explanation of the water cycle algorithm. Let us assume the $1 \times D$ dimensional array as a candidate solution (i.e., stream)
\begin{equation}
\label{FIWCA:Candidate_Sol}
X=[x_1,x_2,x_3,...,x_D],
\end{equation}
where $D$ is the number of design variable. Hence, the initial randomly generated population can be represented by an $N_{pop} \times D$ matrix as given follows
\begin{equation}
\label{FIWCA:Total_population}
\begin{split}
&\text{Total Population}=
\begin{bmatrix}
Sea\\
\ {River}_1 \\
\ {River}_2 \\
\ {River}_3 \\
\ \vdots \\
\ {Stream}_{N_{sr}+1} \\
\ {Stream}_{N_{sr}+1} \\
\ {Stream}_{N_{sr}+3} \\
\ \vdots \\
\ {Stream}_{N_{pop}} \\
\end{bmatrix}\\
&\quad \quad \quad \quad=
\begin{bmatrix}
 x_{1}^{1} & x_{2}^{1} & x_{3}^{1} & \ldots & x_{D}^{1}  \\
\ x_{1}^{2} & x_{2}^{2} & x_{3}^{2} & \ldots & x_{D}^{1}  \\
\ \vdots & \vdots & \vdots & \vdots & \vdots  \\
\ x_{1}^{N_{pop}} & x_{2}^{N_{pop}} & x_{3}^{N_{pop}} & \ldots & x_{D}^{N_{pop}}  \\
\end{bmatrix},\\
\end{split}
\end{equation}
where $N_{pop}$ is the population size. To choose the rivers and sea, first of all the cost value of every generated population, i.e. every stream, need to be determined. The cost of every stream $X^i$ is calculated by evaluating the problem cost function for that stream as given follows
\begin{equation}
\label{FIWCA:Cost}
\begin{aligned}
& C_i={Cost}_i=f(x_1^i,x_2^i,\ldots,x_D^i) & i=1,2,\ldots,N_{pop}.
\end{aligned}
\end{equation}
A user parameter which is summation of number of rivers and sea so called $N_{sr}$ is selected. Among this set, the individual stream with the smallest objective value is selected as the sea and the others are chosen as the rivers, given in the following equations
\begin{equation}
\label{FIWCA:Nsr}
N_{sr}=\text{Number of Rivers}+\underbrace{1}_{Sea},
\end{equation}
\begin{equation}
\label{FIWCA:Nstream}
N_{Stream}=N_{pop}-N_{sr}.
\end{equation}
Therefore, the population of streams that flow towards the rivers and sea is given as follows
\begin{equation}
\footnotesize
\label{FIWCA:Population_of_Streams}
\begin{split}
&\text{ Population of Streams}=
\begin{bmatrix}
 {Stream}_1 \\
\ {Stream}_2 \\
\ \vdots \\
\ {Stream}_{N_{Stream}} \\
\end{bmatrix}\\
&\quad\quad\quad\quad=
\begin{bmatrix}
 x_{1}^{1} & x_{2}^{1} &  \ldots & x_{D}^{1}  \\
\ x_{1}^{2} & x_{2}^{2} &  \ldots & x_{D}^{1}  \\
\ \vdots & \vdots & \vdots  & \vdots  \\
\ x_{1}^{N_{Stream}} & x_{2}^{N_{Stream}} &  \ldots & x_{D}^{N_{Stream}}  \\
\end{bmatrix}.\\
\end{split}
\normalsize
\end{equation}

In fact \eqref{FIWCA:Population_of_Streams} which is a part of the total population with the size of $N_{Stream}$, \eqref{FIWCA:Nstream}. The sea and every river absorb water from the streams. As a matter of fact, water is flowing from streams toward the rivers and the sea by its own nature. The amount of water flow that enters a river or/and sea is different from one stream to another. The rivers also flow towards the sea which is in the deepest location. Therefore, the group of streams designated for every river and sea is determined by the following equation \cite{sadollah2015water}
\begin{equation}
\label{FIWCA:C_n}
\begin{aligned}
& C_n = {Cost}_n - {Cost}_{N_{sr}+1} & n=1,2,\ldots,N_{sr},
\end{aligned}
\end{equation}
\begin{equation}
\label{FIWCA:NS_n}
\begin{aligned}
&{NS}_n = Round \Bigg\{ \Big|\frac{C_n}{\sum_{n=1}^{N_{sr}} C_n} \Big| \times N_{Stream} \Bigg\} & n=1,\ldots,N_{sr},
\end{aligned}
\end{equation}
where ${NS}_n$ is the number of streams which move toward the specific rivers and sea, i.e. ${NS}_1$ is the number of streams move toward the sea and ${NS}_2,\ldots,{NS}_{N_{sr}}$ are the number of streams that flow toward their corresponding rivers.

\begin{figure}
 \centering
 \includegraphics[width=8cm, height=12cm,keepaspectratio]{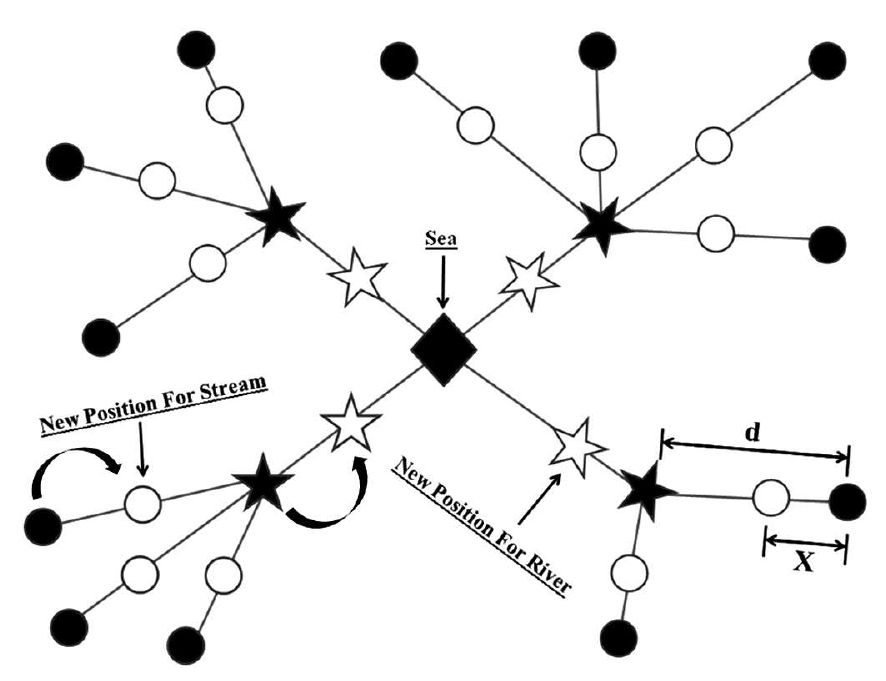}
 \caption{Schematic View of the WCA Optimization Process \cite{eskandar2012water}.}
 \label{Schematic_of_WCA}
 \end{figure}

As it can be observed in nature, raining and precipitation build the streams and the rivers are created by joining the streams to each other. Some of the streams might even move to the sea, directly. Eventually, all the streams and rivers reach the sea which is considered as the current best solution. Fig.~\ref{Schematic_of_WCA} illustrates the process of WCA in which diamond, stars and circles  represent the sea, rivers and streams, respectively. The black circles and stars show the current position of streams and rivers and the white ones denote their new position at the next iteration \cite{eskandar2012water}.

\subsection{Fully Informed Water Cycle Algorithm}

\begin{figure}
\centering
\centering
\subfloat[New Position of the Streams which Directly Flow toward the Sea.]{\includegraphics[totalheight=0.19\textheight]{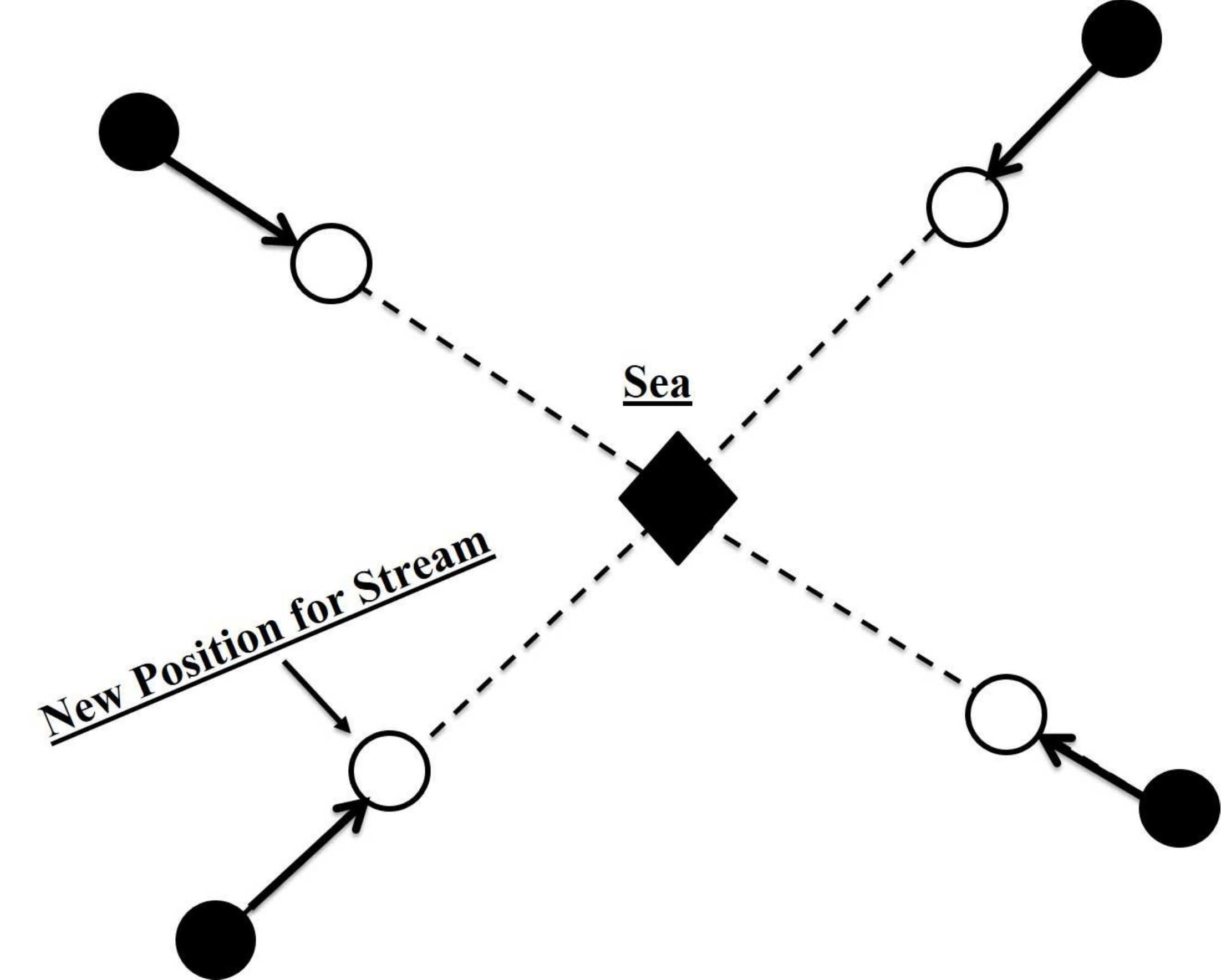}\label{FIWCA_A}}\hspace{25mm}\\
\subfloat[New Position of the Streams which Flow toward the Rivers.]{\includegraphics[totalheight=0.19\textheight]{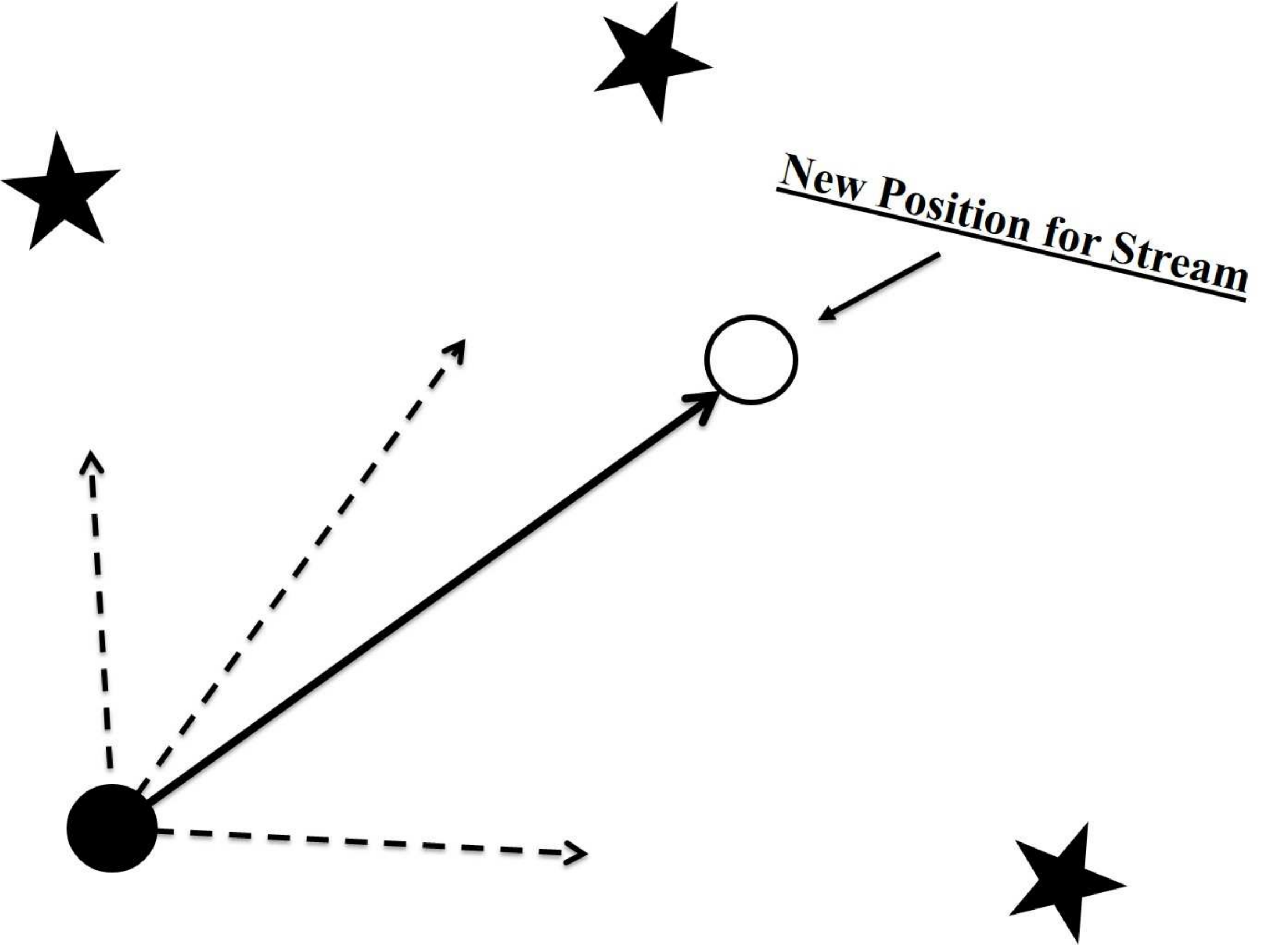}\label{FIWCA_B}}\hspace{25mm}\\
\subfloat[New Position of Every River by Considering the Influence of Sea and all other Rivers.]{\includegraphics[totalheight=0.19\textheight]{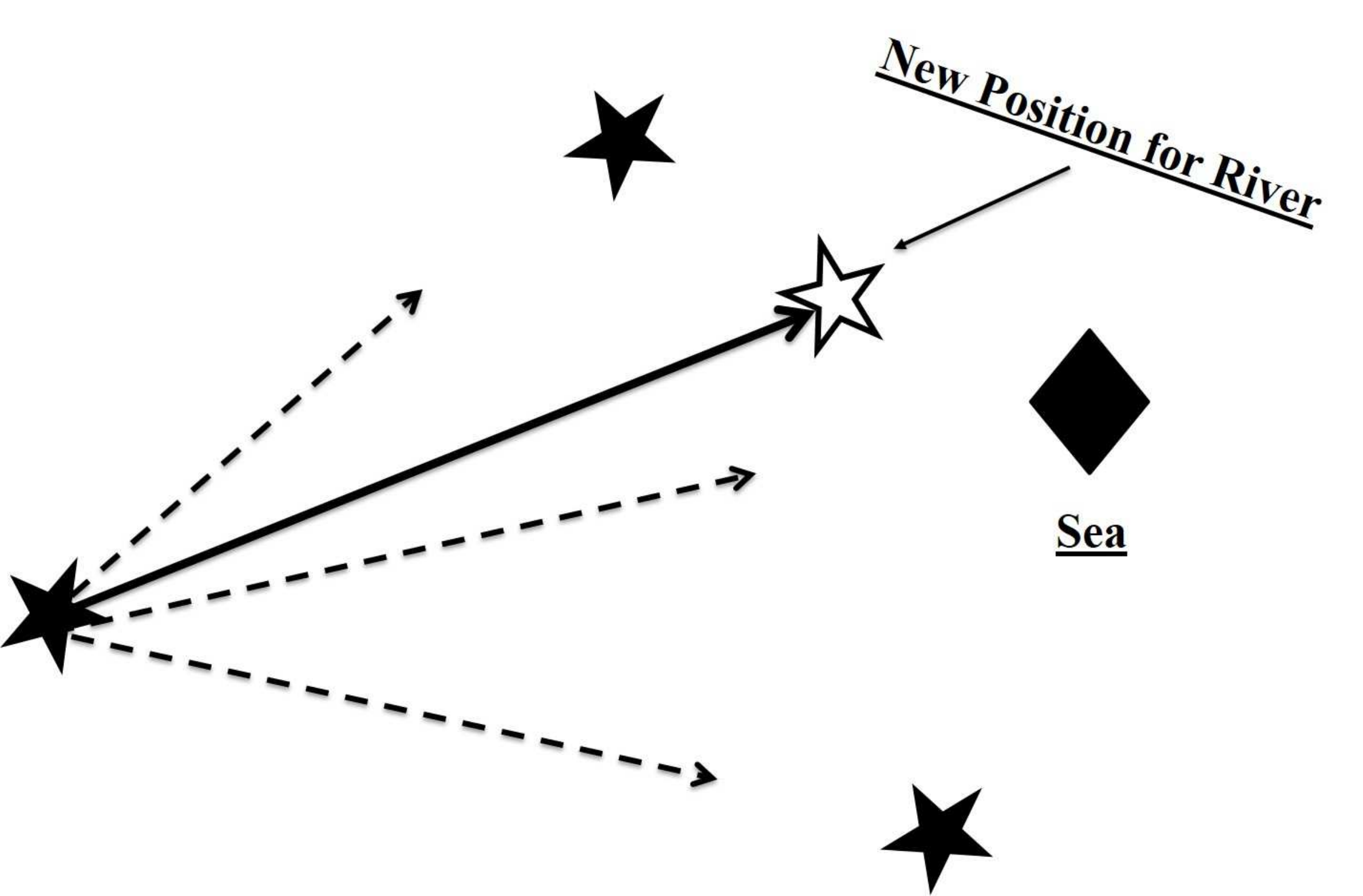}\label{FIWCA_C}}\hspace{25mm}\\
\caption{Schematic View of the Fully Informed Water Cycle Algorithm (FIWCA) Optimization Process.}
\label{Schematic_of_FIWCA}
\end{figure}

In the standard WCA, each stream corresponds to one of the rivers and/or sea without considering the influence of other rivers and the sea. The rivers also just consider the sea as their reference to exchange information for flowing. In this section, a fully informed water cycle algorithm has been proposed which every stream receives the influence of all rivers and/or sea as the information reference for the flowing and updating the new positions. The information in rivers flows not only to the sea, but these information flows can also be affected by all rivers and the sea to update their new positions and to intensively exploit the searching space. Fig.~\ref{Schematic_of_FIWCA} displays the schematic view of the movement process used in FIWCA.

Now, mathematical explanations of updating the streams and rivers positions can be expressed. Let us assume $N_{pop}$ is the total number of streams, of which $N_{sr}-1$ are selected as the rivers and one individual is chosen as the sea. Let $\overrightarrow{(\cdot)}$ represent a vector, the new position of streams and rivers are then given as follows

\begin{subequations}
\small
\label{FIWCA:positions}
\begin{align}
\label{FIWCA:positions_stream_sea}
\begin{split}
& \overrightarrow{X}_{Stream}^i(t+1) = \overrightarrow{X}_{Stream}^i(t) \\
& \quad \quad  + rand \times C \times \big( \overrightarrow{X}_{Sea}(t) - \overrightarrow{X}_{Stream}^{i}(t) \big) \\
& \qquad \qquad \qquad \qquad \qquad \qquad  i=1,2,\ldots,NS_{1},\\
\end{split}\\
\label{FIWCA:positions_stream_river}
\begin{split}
& \overrightarrow{X}_{Stream}^i(t+1) = \overrightarrow{X}_{Stream}^i(t) \\
& \quad \quad    + \sum_{\substack{j=1 \\ j\neq i}}^{N_{sr}-1}rand \times C \times \big( \overrightarrow{X}_{River}^{j}(t) - \overrightarrow{X}_{Stream}^{i}(t) \big)\\
& \qquad \qquad \qquad \qquad \qquad \qquad  i=1,2,\ldots,NS_{2-N_{sr}},\\
\end{split}\\
\label{FIWCA:positions_river_sea}
\begin{split}
& \overrightarrow{X}_{River}^i(t+1) = \overrightarrow{X}_{River}^i(t)\\
& \quad \quad  + rand \times C \times \big( \overrightarrow{X}_{Sea}(t) - \overrightarrow{X}_{River}^{i}(t) \big)\\
& \quad \quad   + \sum_{\substack{j=1 \\ j\neq i}}^{N_{sr}-1}rand \times C \times \big( \overrightarrow{X}_{River}^{j}(t) - \overrightarrow{X}_{River}^{i}(t) \big)\\
& \qquad \qquad \qquad \qquad \qquad \qquad  i=1,2,\ldots,(N_{sr}-1),\\
\end{split}
\end{align}
\normalsize
\end{subequations}
where $1<C<2$ as proposed in \cite{eskandar2012water}. For the best performance $C$ can be equal to $2$,``$rand$'' is a uniformly distributed random number between $0$ and $1$, $t$ is the iteration index, and $i$ is the number of individuals. Note that \eqref{FIWCA:positions_stream_sea} defines the new position for the streams that directly move toward the sea, see Fig.~\ref{FIWCA_A}, in the same way as the one used in standard WCA \cite{eskandar2012water} and \eqref{FIWCA:positions_stream_river} shows the new position for the streams which flow toward the rivers, not only its corresponding river, see Fig.~\ref{FIWCA_B}. In fact, in \eqref{FIWCA:positions_stream_river}, a stream flowing to the river, has been influenced not only by its specific river, but other rivers are also affecting the stream to provide full information of all the rivers for its new position. \eqref{FIWCA:positions_river_sea} represents the new position of every river by considering the influence of the sea and all other rivers, see Fig.~\ref{FIWCA_C}.

Note that whenever the cost value of any stream becomes better (smaller) than the cost value of its corresponding river, the status of that stream and river are exchanged, i.e. the stream turns into a river and the river becomes a stream. The same replacement will occur between a river and the sea.

In order to avoid getting trapped in local optima, the evaporation process has been introduced \cite{sadollah2015water}. In nature, the evaporation process occurs as the streams and rivers run toward the sea, which leads to new precipitations.
\begin{algorithm}
    \caption{Fully Informed Water Cycle Algorithm}
    \label{alg_FIWCA}
    \begin{algorithmic}[1]
        \renewcommand{\algorithmicrequire}{\textbf{Input:}}
        \renewcommand{\algorithmicensure}{\textbf{Output:}}
        \REQUIRE Set user parameters of the WCA: $N_{pop}$, $N_{sr}$, $d_{max}$, and $\text{Max Iteration}$.
        \ENSURE Optimized objective value and design variables.
        \STATE \textit{Initialisation}: Create a random initial population of streams between $LB$ and $UB$ \eqref{FIWCA:Total_population}.
        \WHILE {``($t \leq Max\ Iteration)''$ or (any assumed stopping condition)}
        \FOR {``$i = 1: N_{pop}$''}
        \STATE Update the position of stream directly toward the sea using \eqref{FIWCA:positions_stream_sea}.
        \STATE Calculate the objective function of the newly generated stream.
        \IF {``$\text{Objective} (New\ Stream) < \text{Objective} (Sea)$''}
        \STATE ``Sea = New Stream''.
        \ENDIF
        \STATE Update the position of stream to its corresponding river while considering the information flow from the other rivers using \eqref{FIWCA:positions_stream_river}.
        \STATE Compute the objective function value of the newly generated stream.
        \IF {``$\text{Objective} (New\ Stream) < \text{Objective} (River)$''}
        \STATE River = New Stream.
        \IF {``$\text{Objective} (River) < \text{Objective} (Sea)$''}
        \STATE ``Sea = River''.
        \ENDIF
        \ENDIF
        \STATE Update the position of river to the sea while considering the information flow from other rivers using \eqref{FIWCA:positions_river_sea}.
        \STATE Compute the objective function value of the newly generated river.
        \IF {``$\text{Objective} (New\ River) < \text{Objective} (Sea)$''}
        \STATE ``Sea = New River''.
        \ENDIF
        \ENDFOR
        \STATE \%\% Evaporation Condition for rivers
        \FOR {``$i = 1: N_{sr}-1$''}
        \IF {``$|\overrightarrow{X}_{Sea}-\overrightarrow{X}_{River}| < d_{max} \text{ or } (rand < 0.1)$''}
        \STATE Create new streams and a river (the best one) using \eqref{FIWCA:Evaporation}.
        \ENDIF
        \ENDFOR
        \STATE \%\% Evaporation Condition for streams who directly flow to the sea
        \FOR {``$i = 1: NS_1$''}
        \IF{``$|\overrightarrow{X}_{Sea}-\overrightarrow{X}_{Stream}| < d_{max}$''}
        \STATE Create New streams using \eqref{FIWCA:Evaporation} and the best one is considered as the sea.
        \ENDIF
        \ENDFOR
        \STATE Reduce the $d_{max}$ using \eqref{FIWCA:d_max}.
        \ENDWHILE
        \RETURN Post-process the results and visualization.
    \end{algorithmic}
\end{algorithm}

The evaporation condition in the WCA is defined by whether any stream or river gets close enough (using Euclidian distance) to the sea. This idea is implemented in the standard WCA and FIWCA by the following criteria for the rivers flowing into the sea as well as the streams flowing into the sea, respectively, given as follows \cite{sadollah2015water}

    \begin{algorithmic}[]
        \IF {$|\overrightarrow{X}_{Sea}-\overrightarrow{X}_{River}| < d_{max}$ or\\
         \quad  $rand < 0.1$, $i=1,2,\ldots,N_{sr}-1$}
        \STATE Perform precipitation process using \eqref{FIWCA:Evaporation}.
        \ENDIF
        \IF {$|\overrightarrow{X}_{Sea}-\overrightarrow{X}_{Stream}| < d_{max}$, $i=1,2,\ldots,NS_1$}
        \STATE Perform precipitation process using \eqref{FIWCA:Evaporation}.
        \ENDIF
    \end{algorithmic}
which states that if any river (stream) gets closer than $d_{max}$ to the sea, then the evaporation process occurs. After the evaporation condition is satisfied, it is required to randomly form the new streams in a different position between lower and upper bounds. This phenomenon happens by the precipitation process which provides a new sub-population for which the best stream becomes the river (sea) and the other streams moves toward this river (sea). Indeed, the evaporation and raining process help in the exploration phase of the WCA. The location of newly generated sub-population of the corresponding evaporated river (stream) can be specified as follows
\begin{equation}
\label{FIWCA:Evaporation}
\overrightarrow{X}_{Stream}(t+1) = \overrightarrow{LB} + rand \times (\overrightarrow{UB} - \overrightarrow{LB}),
\end{equation}
in which $\overrightarrow{LB}$ and $\overrightarrow{UB}$ represent the lower and upper bounds of the optimization problem, respectively. Since large values of $d_{max}$ cause extra exploration and very small values lead to the search gravity around the sea, i.e. more exploitation, thus, this parameter, $d_{max}$, may directly affect the search intensity near the sea. For the sake of convergence, it is recommended in \cite{sadollah2015water} to decrease the value of $d_{max}$ adaptively as given follows
\begin{equation}
\label{FIWCA:d_max}
\begin{split}
&d_{max}(t+1) = d_{max}(t) - \frac{d_{max}(t)}{\text{Max Iteration}} \\
&\qquad \qquad\qquad\qquad\qquad\qquad t=1,2,\ldots,\text{Max Iteration}.
\end{split}
\end{equation}
Furthermore, Algorithm \ref{alg_FIWCA} shows the pseudo-code of FIWCA in details.

\subsection {Constraint Handling and Feasibility Rules}

The OPF problem includes a set of constraints on the independent or design variables, such as generators active power outputs and generators voltages, and a set of constraints on the dependent variables, such as line power flows, generators reactive power outputs, and voltages of non-generators buses. The optimization tools and techniques have to be able to handle all these constraints.

Here, the proposed FIWCA utilizes different procedures to maintain each of these groups of constrains to maintain the feasibility of the candidate solutions at each iteration.

To handle the constraints on the design variables which are directly considered by the FIWCA, the following rules are defined \cite{mezura2008empirical}:
\begin{itemize}
  \item \textbf{Rule 1:} Between any two feasible candidate solutions the one with the best objective function value is preferred.
  \item \textbf{Rule 2:} Any feasible candidate solution is preferred to any infeasible candidate solution.
  \item \textbf{Rule 3:} Infeasible solutions containing slight violation of the constraints (from 0.01 in the first iteration to 0.001 in the last iteration) are considered as feasible candidate solutions.
  \item \textbf{Rule 4:} Between two infeasible solution the one with the lowest sum of constraint violation is preferred.
\end{itemize}

As discussed, the OPF problem also includes a set of constraints on the independent variables which are not directly considered in the proposed algorithm. To overcome these constraints, the penalty function approach is utilized, i.e. a penalty function of the each constraint violation is defined and added to the original objective function with a large penalty multiplier.

\section{Optimization Results}

The proposed FIWCA and WCA have been implemented in OPF problems for the IEEE standard 30-bus and 57-bus test systems to obtain better optimal solutions using the proposed improved version and also for the sake of comparison among previous works. The optimization methods considered in this paper have been coded in ``MATLAB'' and solved on an ``Intel 3.40 GHz'' computer with ``8 GB of RAM''. For the original WCA and the fully informed WCA, the population size and the maximum iteration number are fixed to 200 and 100, respectively, i.e. the maximum number of function evaluation of 20,000, in order to have a fair comparison with other algorithms. All numerical test systems are performed for 50 independent replications. Note that in the FIWCA and WCA, the total number of sea and rivers, $N_{sr}$, is considered as 5.

In the literature, a variety of meta-heuristic optimization techniques have been studied for the OPF problems. In this work, the following optimizers have been considered: shuffled frog-leaping algorithm (SFLA) \cite{ghasemi2015improved}, constant inertia particle swarm optimization (CI-PSO) \cite{khorsandi2013modified}, dragonfly algorithm and particle swarm optimization (DA-PSO) \cite{khunkitti2018hybrid}, teaching-learning based optimization (TLBO) and Levy TLBO (LTLBO) \cite{ghasemi2015improved}, multi-objective TLBO (MOTLBO) \cite{nayak2012application},  particle swarm optimization and gravitational search algorithm (PSOGSA) \cite{radosavljevic2015optimal}, grey wolf optimizer (GWO) and developed grey wolf optimizer (DGWO) \cite{abdo2018solving}, enhanced differential evolution (Enhanced DE) \cite{el2015single}, artificial bee colony (ABC) \cite{adaryani2013artificial}, Jaya algorithm \cite{warid2016optimal}, evolving ant direction differential evolution (EADDE), efficient evolutionary algorithm (EEA), enhanced GA for decoupled quadratic load flow (EGA-DQLF) \cite{kumari2010enhanced}, glowworm swarm optimization (GSO) \cite{reddy2016optimal}, hybrid differential evolution and harmony search algorithm (Hybrid DE-HS) \cite{reddy2018optimal},  modified sine-cosine algorithm(MSCA) \cite{attia2018optimal}, bare-bones multi-objective particle swarm optimization (BB-MOPSO) \cite{ghasemi2014multi},  biogeography-based optimization (BBO) \cite{bhattacharya2011application}, chaotic self-adaptive differential harmony search algorithm (CSA-DHS) \cite{arul2013solving}, modified non-dominated sorting genetic algorithm (MNSGA-II) \cite{ghasemi2014multi}, hybrid differential evolution and pattern search (DE-PS) \cite{gitizadeh2014using}, modified particle swarm optimization (MPSO) \cite{karami2015fuzzy,ben2018hybrid}, improved particle swarm optimization (IPSO) \cite{niknam2012improved}, differential search algorithm (DSA) \cite{abaci2016differential}, hybrid modified PSO-SFLA (HMPSO-SFLA) \cite{narimani2013novel}, gravitational search algorithm (GSA) \cite{duman2012optimal}, linearly decreasing inertia weight PSO (LDI-PSO) \cite{adaryani2013artificial}, harmony search algorithm (HSA) \cite{pandiarajan2016fuzzy}.

\subsection{IEEE 30-bus Test System}

\begin{figure}[!t]
 \centering
 \includegraphics[width=8cm, height=9cm,keepaspectratio]{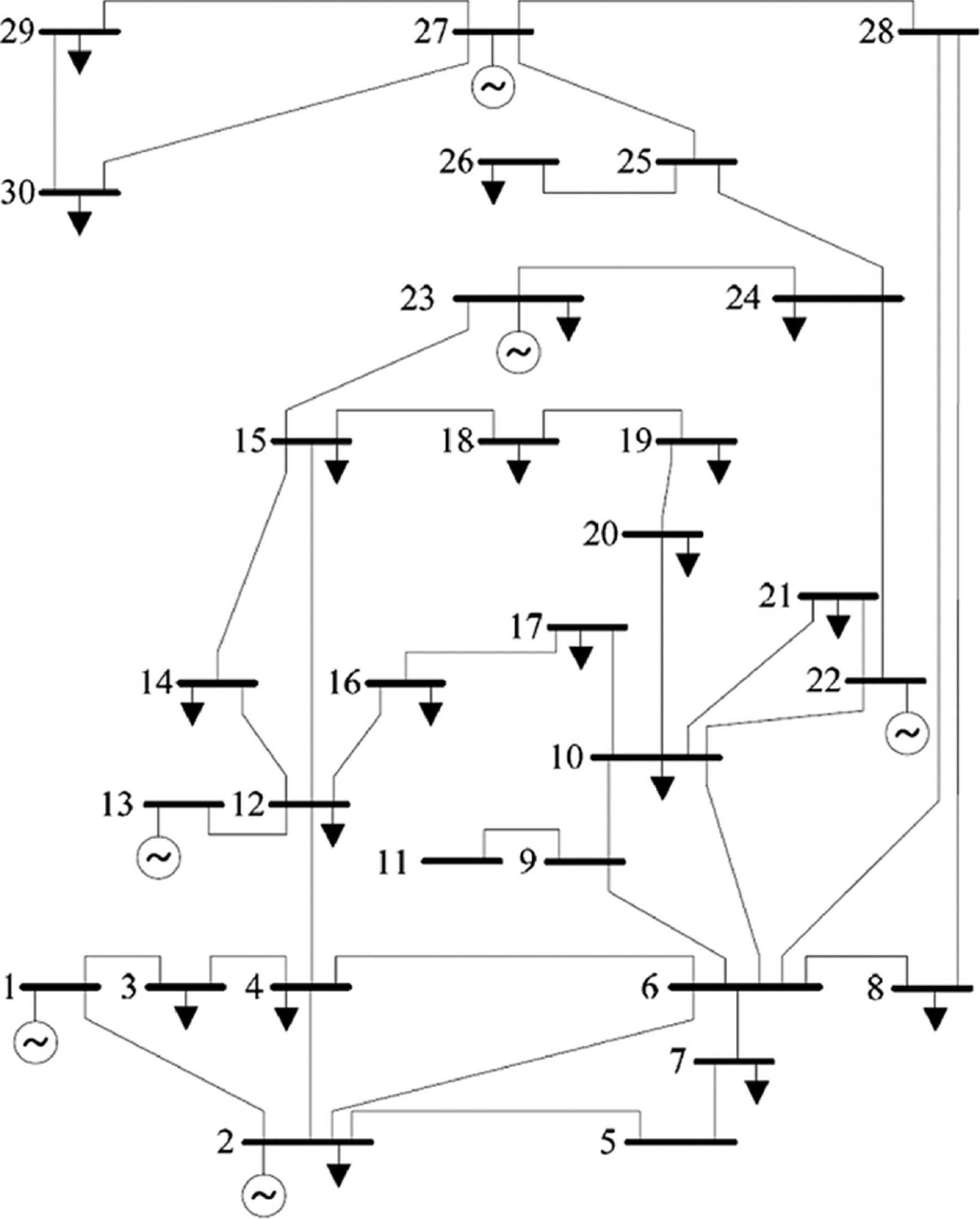}
 \caption{Schematic View of IEEE 30-bus Test System.}
 \label{Schematic_of_30_bus}
 \end{figure}

The IEEE 30-bus test system, as shown in Fig.~\ref{Schematic_of_30_bus}, has six generators fixed at buses 1, 2, 5, 8, 11, and 13, four transformers with the off-nominal tap ratio at lines 6-9, 6-10, 4-12, and 28-27. Load data and line data for the IEEE 30-bus test system can be found in \cite{abido2002optimal}. According to the literature \cite{abido2002optimal}, buses 10, 12, 15, 17, 20, 21, 23, 24, and 29 are considered as shunt VAR compensation buses. The minimum and maximum limits on control variables can be found in Appendix, Table \ref{Case30_Upper_and_Lower_limits}. The proposed FIWCA is implemented on OPF problem with different objective function scenarios to evaluate its effectiveness, as presented in the following subsections.

\begin{sidewaystable*}[]
\centering
\caption{Optimal Control Settings for Different on IEEE standard Cases 30-bus Test System Using the WCA and FIWCA. "N/A" stands for not available.}\label{All Cases_optimal_Values}
{\setlength\arrayrulewidth{1pt}
\begin{tabular}{lllllllllll}
\hline
 \multirow{2}{*}{\textbf{Design Variables}}& \multicolumn{2}{c}{\textbf{Case 1}} & \multicolumn{2}{c}{\textbf{Case 2}}& \multicolumn{2}{c}{\textbf{Case 3}}& \multicolumn{2}{c}{\textbf{Case 4}}& \multicolumn{2}{c}{\textbf{Case 5}}\\
 & \multicolumn{1}{c}{\textbf{WCA}} & \multicolumn{1}{c}{\textbf{FIWCA}}& \multicolumn{1}{c}{\textbf{WCA}} & \multicolumn{1}{c}{\textbf{FIWCA}}& \multicolumn{1}{c}{\textbf{WCA}} & \multicolumn{1}{c}{\textbf{FIWCA}}& \multicolumn{1}{c}{\textbf{WCA}} & \multicolumn{1}{c}{\textbf{FIWCA}} & \multicolumn{1}{c}{\textbf{WCA}} & \multicolumn{1}{c}{\textbf{FIWCA}}\\ \hline
$P_{G1}$ (MW) & $177.0771$ & $177.0756$ & $177.1920$ & $176.0293$ & $139.9999$ & $139.9995$ & $199.5996$ & $199.5996$ & $\hphantom{0}64.9297$ & $\hphantom{0}64.9350$\\
$P_{G2}$ (MW) &	$\hphantom{0}48.6765$	& $\hphantom{0}48.6800$ &	$\hphantom{0}50.0000$	& $\hphantom{0}50.0007$ & $\hphantom{0}54.9998$ & $\hphantom{0}54.9988$ & $\hphantom{0}20.0000$ & $\hphantom20.0000$ & $\hphantom{0}66.4293$ & $\hphantom{0}66.4241$ \\
$P_{G5}$ (MW) &	$\hphantom{0}21.2971$	& $\hphantom{0}21.2965$ &	$\hphantom{0}21.7403$	& $\hphantom{0}21.6417$ & $\hphantom{0}24.1771$ & $\hphantom{0}25.0704$ & $\hphantom{0}21.9998$ & $\hphantom{0}22.5907$ & $\hphantom{0}50.0000$ & $\hphantom{0}50.0000$\\
$P_{G8}$ (MW) &	$\hphantom{0}21.0805$	& $\hphantom{0}21.0806$ &	$\hphantom{0}23.2284$	& $\hphantom{0}22.4090$ & $\hphantom{0}34.9303$ & $\hphantom{0}34.8609$ & $\hphantom{0}26.1736$ & $\hphantom{0}24.5433$ & $\hphantom{0}35.0000$ & $\hphantom{0}35.0000$\\
$P_{G11}$ (MW) &	$\hphantom{0}11.8404$	& $\hphantom{0}11.8390$  &	$\hphantom{0}10.0009$	& $\hphantom{0}12.2134$  & $\hphantom{0}18.8638$ & $\hphantom{0}16.3878$ & $\hphantom{0}12.7518$ & $\hphantom{0}13.5448$ & $\hphantom{0}30.0000$ & $\hphantom{0}30.0000$\\
$P_{G13}$ (MW) &	$\hphantom{0}12.0000$	& $\hphantom{0}12.0000$  &	$\hphantom{0}11.2310$	& $\hphantom{0}11.0339$ & $\hphantom{0}17.4422$ & $\hphantom{0}18.8228$ & $\hphantom{0}12.0616$ & $\hphantom{0}12.2930$ & $\hphantom{0}40.0000$ & $\hphantom{0}40.0000$\\
$V_{G1}$ (p.u.) &	$\hphantom{0}\hphantom{0}1.1000$	& $\hphantom{0}\hphantom{0}1.0999$  &	 $\hphantom{0}\hphantom{0}1.0347$	& $\hphantom{0}\hphantom{0}1.0369$ & $\hphantom{0}\hphantom{0}1.0500$ & $\hphantom{0}\hphantom{0}1.0862$ & $\hphantom{0}\hphantom{0}1.1000$ & $\hphantom{0}\hphantom{0}1.1000$ & $\hphantom{0}\hphantom{0}1.1000$ & $\hphantom{0}\hphantom{0}1.1000$ \\
$V_{G2}$ (p.u.) &	$\hphantom{0}\hphantom{0}1.0877$	& $\hphantom{0}\hphantom{0}1.0877$ &	 $\hphantom{0}\hphantom{0}1.0189$	& $\hphantom{0}\hphantom{0}1.0218$ & $\hphantom{0}\hphantom{0}1.0404$ & $\hphantom{0}\hphantom{0}1.0711$ & $\hphantom{0}\hphantom{0}1.0845$ & $\hphantom{0}\hphantom{0}1.0846$& $\hphantom{0}\hphantom{0}1.0956$ & $\hphantom{0}\hphantom{0}1.0957$\\
$V_{G5}$ (p.u.) &	$\hphantom{0}\hphantom{0}1.0614$	& $\hphantom{0}\hphantom{0}1.0614$ &	 $\hphantom{0}\hphantom{0}1.0105$	& $\hphantom{0}\hphantom{0}1.0141$ & $\hphantom{0}\hphantom{0}1.0157$ & $\hphantom{0}\hphantom{0}1.0377$ & $\hphantom{0}\hphantom{0}1.0589$ & $\hphantom{0}\hphantom{0}1.0592$ & $\hphantom{0}\hphantom{0}1.0781$ & $\hphantom{0}\hphantom{0}1.0781$\\
$V_{G8}$ (p.u.) &	$\hphantom{0}\hphantom{0}1.0693$	& $\hphantom{0}\hphantom{0}1.0693$ &	 $\hphantom{0}\hphantom{0}1.0056$	& $\hphantom{0}\hphantom{0}1.0071$ & $\hphantom{0}\hphantom{0}1.0260$ & $\hphantom{0}\hphantom{0}1.0536$ & $\hphantom{0}\hphantom{0}1.0684$ & $\hphantom{0}\hphantom{0}1.0682$ & $\hphantom{0}\hphantom{0}1.0855$ & $\hphantom{0}\hphantom{0}1.0856$\\
$V_{G11}$ (p.u.) &	$\hphantom{0}\hphantom{0}1.1000$	& $\hphantom{0}\hphantom{0}1.1000$  &	 $\hphantom{0}\hphantom{0}1.0273$	& $\hphantom{0}\hphantom{0}0.9973$ & $\hphantom{0}\hphantom{0}1.1000$ & $\hphantom{0}\hphantom{0}1.0152$ & $\hphantom{0}\hphantom{0}1.1000$ & $\hphantom{0}\hphantom{0}1.1000$ & $\hphantom{0}\hphantom{0}1.1000$ & $\hphantom{0}\hphantom{0}1.1000$\\
$V_{G13}$ (p.u.) &	$\hphantom{0}\hphantom{0}1.0999$	& $\hphantom{0}\hphantom{0}1.0999$  &	 $\hphantom{0}\hphantom{0}0.9978$	& $\hphantom{0}\hphantom{0}0.9977$ & $\hphantom{0}\hphantom{0}1.0909$ & $\hphantom{0}\hphantom{0}1.0506$ & $\hphantom{0}\hphantom{0}1.0999$ & $\hphantom{0}\hphantom{0}1.1000$ & $\hphantom{0}\hphantom{0}1.1000$ & $\hphantom{0}\hphantom{0}1.1000$\\
$T_{6-9}$ (p.u.) &	$\hphantom{0}\hphantom{0}1.0311$	& $\hphantom{0}\hphantom{0}1.0311$  &	 $\hphantom{0}\hphantom{0}1.0433$	& $\hphantom{0}\hphantom{0}1.0111$ & $\hphantom{0}\hphantom{0}0.9247$ & $\hphantom{0}\hphantom{0}0.9949$ & $\hphantom{0}\hphantom{0}1.0297$ & $\hphantom{0}\hphantom{0}1.0297$& $\hphantom{0}\hphantom{0}1.0343$ & $\hphantom{0}\hphantom{0}1.0340$\\
$T_{6-10}$ (p.u.) &	$\hphantom{0}\hphantom{0}0.9000$	& $\hphantom{0}\hphantom{0}0.9000$ &	 $\hphantom{0}\hphantom{0}0.9000$	& $\hphantom{0}\hphantom{0}0.9000$ & $\hphantom{0}\hphantom{0}1.0205$ & $\hphantom{0}\hphantom{0}1.0065$ & $\hphantom{0}\hphantom{0}0.9000$ & $\hphantom{0}\hphantom{0}0.9000$& $\hphantom{0}\hphantom{0}0.9000$ & $\hphantom{0}\hphantom{0}0.9000$\\
$T_{4-12}$ (p.u.) &	$\hphantom{0}\hphantom{0}0.9679$	& $\hphantom{0}\hphantom{0}0.9679$ &	 $\hphantom{0}\hphantom{0}0.9540$	& $\hphantom{0}\hphantom{0}0.9541$ & $\hphantom{0}\hphantom{0}0.9563$ & $\hphantom{0}\hphantom{0}1.0343$ & $\hphantom{0}\hphantom{0}0.9678$ & $\hphantom{0}\hphantom{0}0.9666$& $\hphantom{0}\hphantom{0}0.9597$ & $\hphantom{0}\hphantom{0}0.9600$\\
$T_{28-27}$ (p.u.) &	$\hphantom{0}\hphantom{0}0.9584$	& $\hphantom{0}\hphantom{0}0.9583$ &	 $\hphantom{0}\hphantom{0}0.9694$	& $\hphantom{0}\hphantom{0}0.9688$ & $\hphantom{0}\hphantom{0}0.9281$ & $\hphantom{0}\hphantom{0}0.9997$ & $\hphantom{0}\hphantom{0}0.9572$ & $\hphantom{0}\hphantom{0}0.9570$& $\hphantom{0}\hphantom{0}0.9623$ & $\hphantom{0}\hphantom{0}0.9625$\\
$Q_{C10}$ (MVAR)	& $\hphantom{0}\hphantom{0}5.0000$	& $\hphantom{0}\hphantom{0}5.0000$ &	 $\hphantom{0}\hphantom{0}5.0000$	& $\hphantom{0}\hphantom{0}4.9877$ & \hphantom{0}\hphantom{0}\hphantom{0}N/A	& \hphantom{0}\hphantom{0}\hphantom{0}N/A & $\hphantom{0}\hphantom{0}5.0000$ & $\hphantom{0}\hphantom{0}5.0000$& $\hphantom{0}\hphantom{0}5.0000$ & $\hphantom{0}\hphantom{0}5.0000$\\
$Q_{C12}$ (MVAR)	& $\hphantom{0}\hphantom{0}5.0000$	& $\hphantom{0}\hphantom{0}5.0000$ &	 $\hphantom{0}\hphantom{0}0.0000$	& $\hphantom{0}\hphantom{0}0.0000$ & \hphantom{0}\hphantom{0}\hphantom{0}N/A	& \hphantom{0}\hphantom{0}\hphantom{0}N/A & $\hphantom{0}\hphantom{0}5.0000$ & $\hphantom{0}\hphantom{0}5.0000$ & $\hphantom{0}\hphantom{0}5.0000$ & $\hphantom{0}\hphantom{0}5.0000$\\
$Q_{C15}$ (MVAR)	& $\hphantom{0}\hphantom{0}5.0000$	& $\hphantom{0}\hphantom{0}5.0000$ &	 $\hphantom{0}\hphantom{0}5.0000$	& $\hphantom{0}\hphantom{0}5.0000$ & \hphantom{0}\hphantom{0}\hphantom{0}N/A	& \hphantom{0}\hphantom{0}\hphantom{0}N/A & $\hphantom{0}\hphantom{0}5.0000$ & $\hphantom{0}\hphantom{0}4.9450$ & $\hphantom{0}\hphantom{0}5.0000$ & $\hphantom{0}\hphantom{0}5.0000$ \\
$Q_{C17}$ (MVAR)	& $\hphantom{0}\hphantom{0}5.0000$	& $\hphantom{0}\hphantom{0}5.0000$ &	 $\hphantom{0}\hphantom{0}0.0213$	& $\hphantom{0}\hphantom{0}0.0000$ & \hphantom{0}\hphantom{0}\hphantom{0}N/A	& \hphantom{0}\hphantom{0}\hphantom{0}N/A & $\hphantom{0}\hphantom{0}5.0000$ & $\hphantom{0}\hphantom{0}5.0000$ & $\hphantom{0}\hphantom{0}5.0000$ & $\hphantom{0}\hphantom{0}5.0000$\\
$Q_{C20}$ (MVAR)	& $\hphantom{0}\hphantom{0}4.3034$	& $\hphantom{0}\hphantom{0}4.3069$ &	 $\hphantom{0}\hphantom{0}5.0000$	& $\hphantom{0}\hphantom{0}5.0000$ & \hphantom{0}\hphantom{0}\hphantom{0}N/A	& \hphantom{0}\hphantom{0}\hphantom{0}N/A & $\hphantom{0}\hphantom{0}5.0000$ & $\hphantom{0}\hphantom{0}4.9677$ & $\hphantom{0}\hphantom{0}4.1321$ & $\hphantom{0}\hphantom{0}4.1505$\\
$Q_{C21}$ (MVAR)	& $\hphantom{0}\hphantom{0}5.0000$	& $\hphantom{0}\hphantom{0}5.0000$ &	 $\hphantom{0}\hphantom{0}5.0000$	& $\hphantom{0}\hphantom{0}5.0000$ & \hphantom{0}\hphantom{0}\hphantom{0}N/A	& \hphantom{0}\hphantom{0}\hphantom{0}N/A & $\hphantom{0}\hphantom{0}5.0000$ & $\hphantom{0}\hphantom{0}5.0000$& $\hphantom{0}\hphantom{0}5.0000$ & $\hphantom{0}\hphantom{0}5.0000$ \\
$Q_{C23}$(MVAR)	& $\hphantom{0}\hphantom{0}2.6507$	& $\hphantom{0}\hphantom{0}2.6422$ &	 $\hphantom{0}\hphantom{0}5.0000$	& $\hphantom{0}\hphantom{0}5.0000$ & \hphantom{0}\hphantom{0}\hphantom{0}N/A	& \hphantom{0}\hphantom{0}\hphantom{0}N/A & $\hphantom{0}\hphantom{0}5.0000$ & $\hphantom{0}\hphantom{0}2.5560$& $\hphantom{0}\hphantom{0}2.5322$ & $\hphantom{0}\hphantom{0}2.5603$\\
$Q_{C24}$ (MVAR)	& $\hphantom{0}\hphantom{0}5.0000$	& $\hphantom{0}\hphantom{0}5.0000$ &	 $\hphantom{0}\hphantom{0}5.0000$	& $\hphantom{0}\hphantom{0}5.0000$& \hphantom{0}\hphantom{0}\hphantom{0}N/A	& \hphantom{0}\hphantom{0}\hphantom{0}N/A & $\hphantom{0}\hphantom{0}0.0003$ & $\hphantom{0}\hphantom{0}5.0000$& $\hphantom{0}\hphantom{0}5.0000$ & $\hphantom{0}\hphantom{0}5.0000$\\
$Q_{C29}$ (MVAR)	& $\hphantom{0}\hphantom{0}2.3028$	& $\hphantom{0}\hphantom{0}2.3045$ &	 $\hphantom{0}\hphantom{0}2.6755$	& $\hphantom{0}\hphantom{0}2.6105$ & \hphantom{0}\hphantom{0}\hphantom{0}N/A	& \hphantom{0}\hphantom{0}\hphantom{0}N/A & $\hphantom{0}\hphantom{0}4.9994$ & $\hphantom{0}\hphantom{0}2.2920$& $\hphantom{0}\hphantom{0}2.1391$ & $\hphantom{0}\hphantom{0}2.1450$\\ \hline
\textbf{Fuel cost (\$/h)}	& $798.8608$	& $798.8608$ &	$803.9741$	& $803.9370$ & $647.5310$ & $646.6892$ & $917.0858$ & $917.0740$ & $941.8130$ & $941.8048$\\ \hline
\textbf{Emission (ton/h)} & $\hphantom{0}\hphantom{0}0.3661$	& $\hphantom{0}\hphantom{0}0.3661$ &	 $\hphantom{0}\hphantom{0}0.3671$	& $\hphantom{0}\hphantom{0}0.3634$ & $\hphantom{0}\hphantom{0}0.2834$ & $\hphantom{0}\hphantom{0}0.2835$ & $\hphantom{0}\hphantom{0}0.4405$ & $\hphantom{0}\hphantom{0}0.4401$ & $\hphantom{0}\hphantom{0}0.2047$ & $\hphantom{0}\hphantom{0}0.2047$ \\ \hline
\textbf{Power loss (MW)}	& $\hphantom{0}\hphantom{0}8.5718$	& $\hphantom{0}\hphantom{0}8.5719$ &	 $\hphantom{0}\hphantom{0}9.9928$	& $\hphantom{0}\hphantom{0}9.9282$ & $\hphantom{0}\hphantom{0}7.0133$ & $\hphantom{0}\hphantom{0}6.7404$ & $\hphantom{0}\hphantom{0}9.1866$ & $\hphantom{0}\hphantom{0}9.1715$ & $\hphantom{0}\hphantom{0}2.9591$ & $\hphantom{0}\hphantom{0}2.9591$ \\ \hline
\textbf{VD (p.u.)}	& $\hphantom{0}\hphantom{0}2.0379$	& \hphantom{0}\hphantom{0}$2.0381$ &	 $\hphantom{0}\hphantom{0}0.0974$	& $\hphantom{0}\hphantom{0}0.0969$ & $\hphantom{0}\hphantom{0}1.1870$ & $\hphantom{0}\hphantom{0}0.4675$ & $\hphantom{0}\hphantom{0}1.9943$ & $\hphantom{0}\hphantom{0}2.0398$ & $\hphantom{0}\hphantom{0}2.2653$ & $\hphantom{0}\hphantom{0}2.2671$ \\ \hline
\end{tabular}}
\normalsize
\end{sidewaystable*}

\subsubsection{\textbf{Case 1:} Minimization of Quadratic Fuel Cost}

In this case, the objective function of the OPF problem is to minimize the total generation fuel cost by defining a quadratic objective function $J_1$ as given follows
\begin{equation}
\label{FIWCA:Case1_cost}
J_1 = \sum_{i=1}^{NG} f_i(P_{Gi}) = \sum_{i=1}^{NG} \big(c_{i0} + c_{i1} P_{Gi} + c_{i2} P_{Gi}^2 \big),
\end{equation}
\begin{table}[]
\centering
\caption{Comparison of Statistical Optimization Results for Case 1, i.e. Fuel Costs (\$/h), Using Several Optimizers.}\label{Case1_comparison}
{\setlength\arrayrulewidth{1pt}
\begin{tabular}{llll}
\hline
\textbf{Algorithms} & \multicolumn{1}{l}{\textbf{\begin{tabular}[c]{@{}l@{}}Minimum\\ Fuel Cost\end{tabular}}} & \multicolumn{1}{l}{\textbf{\begin{tabular}[c]{@{}l@{}}Average \\ Fuel Cost\end{tabular}}} & \multicolumn{1}{l}{\textbf{\begin{tabular}[c]{@{}l@{}}Maximum\\ Fuel Cost\end{tabular}}} \\ \hline
\textbf{SFLA \cite{ghasemi2015improved}} & $802.2100$ & N/A & N/A\\
\textbf{DA-PSO \cite{khunkitti2018hybrid}} &	$802.1241$ & N/A	& N/A \\
\textbf{MTLBO \cite{shabanpour2014modified}} &	$801.8925$ & $801.95$	& N/A \\
\textbf{PSO \cite{ghasemi2015improved} }& $801.8900$ & N/A & N/A\\
\textbf{PSOGSA \cite{radosavljevic2015optimal}} &	$801.4986$ & N/A	& N/A \\
\textbf{GWO \cite{abdo2018solving}} &	$801.2590$ & $802.6630$	& $804.8980$ \\
\textbf{Enhanced DE \cite{el2015single}} &	$801.2300$ & N/A	& N/A \\
\textbf{ABC \cite{adaryani2013artificial}} & $800.6600$ & $800.8715$ & $801.8674$\\
\textbf{Jaya \cite{warid2016optimal}} &	$800.4986$ & N/A	& N/A \\
\textbf{DGWO \cite{abdo2018solving}} &	$800.4330$ & $800.4674$	& $800.4989$ \\
\textbf{EADDE \cite{vaisakh2011evolving} }& $800.2041$ & N/A & N/A\\
\textbf{EEA \cite{reddy2014faster} }& $800.0831$ & $800.1730$ & $800.2123$\\
\textbf{EGA-DQLF \cite{kumari2010enhanced}} & $799.5600$ & N/A & N/A\\  \hline
\textbf{FIWCA (This Study)} & $\textbf{798.8608}$ & $\textbf{798.8609}$ & $\textbf{798.8612}$\\ \hline
\end{tabular}}
\normalsize
\end{table}
where $P_{Gi}$ and $f_i$ are the output power and the generation fuel cost of \emph{i}-th generator, respectively, and $c_{i0}$, $c_{i1}$, and $c_{i2}$ are the cost coefficients of the $i$-th generator which can be found in Appendix, Table \ref{Case1:cost_coefficients}. The optimal values of design variables for Case 1, achieved by the proposed FIWCA, can be found in Table \ref{All Cases_optimal_Values}. The total generation fuel cost achieved by FIWCA is $\$798.8608$ per hour, while the total emission and power losses are $0.3661$ ton per hour and $8.5719$ MW, respectively. Comparison of the optimized fuel cost using the FIWCA along with the other reported optimization methods is tabulated in Table \ref{Case1_comparison}.

By observing Table \ref{Case1_comparison}, the minimum total fuel cost obtained by FIWCA shows better stability in terms of statistical results where the average fuel cost for FIWCA is $\$798.8609$ per hour. The best attainable results of Case 1 have been highlighted in bold in Table \ref{Case1_comparison}.

\subsubsection{\textbf{Case 2:} Improvement of Voltage Profiles}

One of the critical variables of any power flow system, closely correlated with the stability and security of the power system, is bus voltages of the grid. Hence, improving the quality of the bus voltage profiles and keeping them close to the reference value of 1.0 p.u. can be considered as one of the desired goals of OPF problem. In this case, the improved algorithm has been implemented to minimize the total generation fuel cost, while forcing the bus voltages close to 1.0 p.u.. Thus, the objective function can be defined as given follows
\begin{equation}
\label{FIWCA:Case2_cost}
J_2 = \sum_{i=1}^{NG} \big(c_{i0} + c_{i1} P_{Gi} + c_{i2} P_{Gi}^2 \big) + w \sum_{i=1}^{NPQ} |V_i - 1.0|.
\end{equation}
\begin{table}[]
\centering
\caption{Comparison of Statistical Optimization Results for Case 2, i.e. Total Voltage Deviation (V.D.) (p.u.), Obtained by Several Optimizers.}\label{Case2_comparison}
{\setlength\arrayrulewidth{1pt}
\begin{tabular}{llll}
\hline
\textbf{Algorithms} & \multicolumn{1}{l}{\textbf{\begin{tabular}[c]{@{}l@{}}Minimum\\ V.D.\end{tabular}}} & \multicolumn{1}{l}{\textbf{\begin{tabular}[c]{@{}l@{}}Average \\ V.D.\end{tabular}}} & \multicolumn{1}{l}{\textbf{\begin{tabular}[c]{@{}l@{}}Maximum\\ V.D.\end{tabular}}} \\ \hline
\textbf{ABC \cite{khorsandi2013modified}} &	$0.1351$ & N/A	& N/A \\
\textbf{CI-PSO \cite{khorsandi2013modified}} &	$0.1350$ & N/A	& N/A \\
\textbf{MABC \cite{khorsandi2013modified}} &	$0.1292$ & N/A	& N/A \\
\textbf{Jaya \cite{warid2016optimal}} &	$0.1243$ & N/A	& N/A \\
\textbf{GWO \cite{mahdad2015blackout}} &	$0.1076$ & N/A	& N/A \\
\textbf{MSCA \cite{attia2018optimal}} &	$0.1031$ & N/A	& N/A \\
\textbf{BB-MOPSO \cite{ghasemi2014multi}} &	$0.1021$ & N/A	& N/A \\
\textbf{BBO \cite{bhattacharya2011application}}	& $0.1020$ &	$0.1105$ &	 $0.1207$ \\
\textbf{CSA-DHS \cite{arul2013solving}} &	$0.1006$ & N/A	& N/A \\
\textbf{MNSGA-II \cite{ghasemi2014multi}} &	$0.0989$ & N/A	& N/A \\  \hline
\textbf{FIWCA (This Study)}	 &	$\textbf{0.0969}$ &	$\textbf{0.0972}$ &	$\textbf{0.0998}$ \\ \hline
\end{tabular}}
\normalsize
\end{table}
The first term of \eqref{FIWCA:Case2_cost} demonstrates the total fuel cost and the second term shows the cost on voltage profile tracking which tries to minimize the total voltage deviation (V.D.) and $w$ is the weighting factor and considered as 100 in this simulation. The FIWCA has been applied to this case to determine the optimal values of design variables by minimizing the fuel cost, as presented in Table \ref{All Cases_optimal_Values}, while keeping the bus voltage profiles close to the desired profile.

The total voltage deviation obtained by the FIWCA is $0.0969$ p.u.. Table \ref{Case2_comparison} demonstrates a comparison of performances among the previous studies and FIWCA on the cost of voltage profiles. The best attainable voltage deviation cost values have been highlighted in bold in Table \ref{Case2_comparison} which shows the minimum total voltage deviation is obtained by FIWCA is smaller than the other algorithms. The obtained statistical results also show better performance of FIWCA with average total voltage deviation of $0.0971$ p.u. which is very reasonable comparing to the other algorithms.

\subsubsection{\textbf{Case 3:} Piecewise Quadratic Fuel Cost Function}

In real power generation units, it is possible that thermal generating units can be supplied by various types of fuel sources, e.g., natural gas, coal, and oil. Therefore, the total fuel cost function of this generation units may be composed of piecewise quadratic cost functions specified by the type of consuming fuel. Here, we suppose that the generation units at buses 1 and 2 may be fed by different fuel types depending on their ranges of working. Hence, the generation fuel cost function of these units can be redefined using the following equation
\begin{equation}
\label{FIWCA:Case3_cost_Gen}
\begin{split}
&f_i(P_{Gi}) =\\
&\begin{cases}
      c_{i0_1} + c_{i1_1} P_{Gi} + c_{i2_1} P_{Gi}^2 &  P_{Gi}^{min}\leq P_{Gi} \leq P_{Gi_1}\\
      c_{i0_2} + c_{i1_2} P_{Gi} + c_{i2_2} P_{Gi}^2 &  P_{Gi_1}\leq P_{Gi} \leq P_{Gi_2}\\
      \quad\quad\quad\quad\vdots  & \quad\quad\quad\quad\vdots\\
      c_{i0_3} + c_{i1_3} P_{Gi} + c_{i2_3} P_{Gi}^2 &  P_{Gi_{k-1}}\leq P_{Gi} \leq P_{Gi}^{max}
\end{cases},\\ &\qquad \qquad\qquad\qquad\qquad\qquad\qquad\qquad\qquad\qquad i=1,2,
\end{split}
\end{equation}
where $c_{i0_k}$, $c_{i1_k}$, and $c_{i2_k}$ are the cost coefficients of $i$-th generator with fuel type $k$ which their values can be seen in Appendix, Table \ref{Case3:cost_coefficients}. Thus, the total generation fuel cost function can be defined as follows
\begin{equation}
\label{FIWCA:Case3_cost}
\begin{split}
&J_3 = \sum_{i=1}^{2} \big(c_{i0_k} + c_{i1_k} P_{Gi} + c_{i2_k} P_{Gi}^2 \big) \\
&\qquad \qquad\qquad + \sum_{i=3}^{NG} \big(c_{i0} + c_{i1} P_{Gi} + c_{i2} P_{Gi}^2 \big).
\end{split}
\end{equation}
\begin{table}[]
\centering
\caption{Comparison of Statistical Optimization Results for Case 3, i.e. Fuel Costs (\$/h), Obtained by Several Optimizers.}\label{Case3_comparison}
{\setlength\arrayrulewidth{1pt}
\begin{tabular}{llll}
\hline
\textbf{Algorithms} & \multicolumn{1}{l}{\textbf{\begin{tabular}[c]{@{}l@{}}Minimum\\ Fuel Cost\end{tabular}}} & \multicolumn{1}{l}{\textbf{\begin{tabular}[c]{@{}l@{}}Average \\ Fuel Cost\end{tabular}}} & \multicolumn{1}{l}{\textbf{\begin{tabular}[c]{@{}l@{}}Maximum\\ Fuel Cost\end{tabular}}} \\ \hline
\textbf{MDE \cite{ben2018hybrid}}	 &	$650.2800$ & 	N/A & 	N/A \\
\textbf{ABC \cite{adaryani2013artificial}}	 &	$649.0855$ & 	$654.0784$ & 	 $659.7708$ \\
\textbf{Fuzzy-GA \cite{hsiao2004optimal}}	 &	$648.2309$ & 	$648.4410$ & 	$648.7682$ \\
\textbf{TLBO \cite{ghasemi2015improved}} &	$647.8125$ & $647.8335$	& $647.8415$ \\
\textbf{BBO \cite{bhattacharya2011application}}	& $647.7437$ & 	$647.7645$ & 	 $647.7928$ \\
\textbf{DE-PS \cite{gitizadeh2014using}} &	$647.5500$ & $647.6029$	& $647.9253$ \\
\textbf{LTLBO \cite{ghasemi2015improved}} &	$647.4315$ & $647.4725$	& $647.8638$ \\  \hline
\textbf{FIWCA (This Study)}	 &	$\textbf{646.6892}$ & $\textbf{647.4275}$ & $\textbf{647.5760}$ \\ \hline
\end{tabular}}
\normalsize
\end{table}

In order to have a fair comparison with the previous works, the upper limit of bus voltages is set to 1.05 and also no VAR compensation bus is considered. As presented in Table \ref{All Cases_optimal_Values}, the total fuel cost obtained by the FIWCA is $\$646.6892$ per hour. Table \ref{Case3_comparison} demonstrates the statistical comparison using different optimizers for solving this specific case of OPF problem, which shows better performance of the FIWCA comparing to the other algorithms in case of the minimum and average cost values over 50 replications, where the average cost value by FIWCA is $\$647.4275$ per hour. The best attainable fuel cost values have been highlighted in bold in Table \ref{Case3_comparison}.

\subsubsection{\textbf{Case 4:} Quadratic Cost function with Valve Point Loading Effect}

To improve the optimal solution of power flow problems with non-convex cost functions, the valve point loading effect on generators can be considered in the generation fuel cost function of the generators. As reported in the literature \cite{abido2002optimal}, the new fuel cost function can be re-written for Generators 1 and 2 and the other generators cost functions would be kept as the same with Case 1.

The new fuel cost function of Generators 1 and 2 can be given as follows

\begin{equation}
\label{FIWCA:Case4_cost_Gen}
\begin{split}
&f_i(P_{Gi}) = c_{i0} + c_{i1} P_{Gi} + c_{i2} P_{Gi}^2 \\
& \qquad\qquad + \Big| d_i \sin \big( e_i ( P_{Gi}^{min} - P_{Gi} ) \big) \Big| \quad\quad i=1,2.
\end{split}
\end{equation}
The cost coefficients of quadratic part of cost function $c_{i0}$, $c_{i1}$ and, $c_{i2}$ and the coefficients of valve point loading effect term, $d_i$ and $e_i$, can be found in Appendix, Table \ref{Case4:cost_coefficients}. So, the modified objective function for Case 4 is defined as
\begin{equation}
\label{FIWCA:Case4_cost}
\begin{split}
&J_4= \\
&  \sum_{i=1}^{2} \Bigg(c_{i0} + c_{i1} P_{Gi} + c_{i2} P_{Gi}^2 + \Big| d_i \sin \big( e_i ( P_{Gi}^{min} - P_{Gi} ) \big) \Big| \Bigg)\\
&\qquad +\sum_{i=3}^{NG} \big(c_{i0} + c_{i1} P_{Gi} + c_{i2} P_{Gi}^2 \big).
\end{split}
\end{equation}

To find the detailed optimized control parameters by WCA and FIWCA for this case readers are referred to Table \ref{All Cases_optimal_Values}. Table \ref{Case4_comparison} tabulates the obtained statistical optimization results using the different optimization algorithms which demonstrate better performance of the WCA and especially FIWCA comparing to the previous algorithms

The best fuel cost value obtained by the FIWCA is $\$917.0740$ per hour. The emission and power loss obtained by the FIWCA are $0.4401$ ton per hour and $9.1715$ MW, respectively, which are lower values comparing to the othe algorithms, i.e. the FIWCA resulted in smaller fuel cost value, emission and power loss. The average cost values obtained by the FIWCA is $\$917.3205$ per hour.

The best-obtained fuel cost values have been highlighted in bold in Table \ref{Case4_comparison}. As can be seen in Table \ref{Case4_comparison}, the only algorithm that could compete in terms of solution quality with the proposed method is BBO with the average of $\$919.8389$ per hour against the value of $\$917.3205$ per hour obtained by the FIWCA. The remaining optimization methods show worse optimal solutions.

\begin{table}[]
\centering
\caption{Comparison of Statistical Optimization Results for Case 4, i.e. Fuel Costs (\$/h), Obtained by Several Optimizers.}\label{Case4_comparison}
{\setlength\arrayrulewidth{1pt}
\begin{tabular}{llll}
\hline
\textbf{Algorithms} & \multicolumn{1}{l}{\textbf{\begin{tabular}[c]{@{}l@{}}Minimum\\ Fuel Cost\end{tabular}}} & \multicolumn{1}{l}{\textbf{\begin{tabular}[c]{@{}l@{}}Average \\ Fuel Cost\end{tabular}}} & \multicolumn{1}{l}{\textbf{\begin{tabular}[c]{@{}l@{}}Maximum\\ Fuel Cost\end{tabular}}} \\ \hline
\textbf{MPSO \cite{karami2015fuzzy,ben2018hybrid}}	 &	$952.3000$ & 	N/A & 	 N/A \\
\textbf{ABC \cite{adaryani2013artificial}}	 &	$945.4495$ & $960.5647$ & $973.5995$ \\
\textbf{MDE \cite{ben2018hybrid}}	 &	$930.94$ & 	N/A & 	N/A \\
\textbf{BBO \cite{bhattacharya2011application}}	& $919.7647$ & $919.8389$ & $919.8876$ \\
\textbf{PSOGSA \cite{radosavljevic2015optimal}} &	$919.65785$ & N/A	& N/A \\
\textbf{TLBO \cite{ghasemi2015improved}} &	$919.3943$ & $919.4710$	& $919.6483$ \\
\textbf{DE-PS \cite{gitizadeh2014using}} &	$919.0175$ & $ 919.4750$	& $919.9842$ \\  \hline
\textbf{FIWCA (This Study)}	 &	$\textbf{917.0740}$ & $\textbf{917.3205}$ & $\textbf{918.1678}$ \\ \hline
\end{tabular}}
\normalsize
\end{table}

\subsubsection{\textbf{Case 5:} Minimization of Fuel Cost and Emission}

In this case, we minimize the emission function for the OPF problem as well as the fuel cost. For this purpose, the emission function of \emph{SOX} and \emph{NOX}, two impressive emitted gas of electrical power generation processes, are considered. The emission function can be presented as given follows
\begin{equation}
\label{FIWCA:Case5_cost_emission}
\begin{split}
&f_{E_i}(P_{Gi}) = c_{i0}^{SOX} + c_{i1}^{SOX} P_{Gi} + c_{i2}^{SOX} P_{Gi}^2 \\
&\qquad\qquad\qquad\qquad + d_{i}^{NOX} \exp (e_{i}^{NOX} P_{Gi}),\\
\end{split}
\end{equation}
\begin{table}[]
\centering
\caption{Comparison of Statistical Optimization Results for Case 5, i.e. Total Emission (ton/h), Obtained by Several Optimizers.}\label{Case5_comparison}
{\setlength\arrayrulewidth{1pt}
\begin{tabular}{llll}
\hline
\textbf{Algorithms} & \multicolumn{1}{l}{\textbf{\begin{tabular}[c]{@{}l@{}}Minimum\\ Emission\end{tabular}}} & \multicolumn{1}{l}{\textbf{\begin{tabular}[c]{@{}l@{}}Average \\ Emission\end{tabular}}} & \multicolumn{1}{l}{\textbf{\begin{tabular}[c]{@{}l@{}}Maximum\\ Emission\end{tabular}}} \\ \hline
\textbf{TLBO \cite{ghasemi2015improved}} &	$0.2124$ & $ 0.2141$	& $0.2156$ \\
\textbf{MDE \cite{ben2018hybrid}}	 &	$0.2093$ & 	N/A & 	N/A \\
\textbf{IPSO \cite{niknam2012improved}} &	$0.2058$ & N/A	& N/A \\
\textbf{DSA \cite{abaci2016differential}} &	$0.2058$ & N/A	& N/A \\
\textbf{DE-PS \cite{gitizadeh2014using}} &	$0.2052$ & $ 0.2071$	& $ 0.2177$ \\
\textbf{HMPSO-SFLA \cite{narimani2013novel}}	 &	$0.2052$	& N/A	& N/A \\
\textbf{DA-PSO \cite{khunkitti2018hybrid}} &	$0.2049$ & N/A	& N/A \\
\textbf{ABC \cite{adaryani2013artificial}}	 &	$0.2048$	& N/A	& N/A \\  \hline
\textbf{FIWCA (This Study)}	 &	$\textbf{0.2047}$ & $\textbf{0.2047}$ & $\textbf{0.2047}$ \\ \hline
\end{tabular}}
\normalsize
\end{table}
where $c_{i0}^{SOX}$, $c_{i1}^{SOX}$, and $c_{i2}^{SOX}$ are the \emph{SOX} emission coefficients, and $d_{i}^{NOX}$ and $e_{i}^{NOX}$ are the \emph{NOX} emission coefficient of \emph{i}-th generator. The values of these coefficients can be found in Appendix, Table \ref{Case5:cost_coefficients}. Hence, the objective function of this case can be expressed as given follows
\begin{equation}
\label{FIWCA:Case4_cost}
\begin{split}
& J_5 = \sum_{i=1}^{NG} \big(c_{i0} + c_{i1} P_{Gi} + c_{i2} P_{Gi}^2 \big) \\
&\qquad + \tau \sum_{i=1}^{NG} \Big( c_{i0}^{SOX} + c_{i1}^{SOX} P_{Gi} + c_{i2}^{SOX} P_{Gi}^2 \\
&\qquad\qquad\qquad + d_{i}^{NOX} \exp (e_{i}^{NOX} P_{Gi})\Big).
\end{split}
\end{equation}
Note that in this simulation, the value of $\tau$ is set to 50000. The optimized value of control variables for this case also can be observed in Table \ref{All Cases_optimal_Values}. Table \ref{Case5_comparison} represents the total amount of emission of generator by optimizing the power flow problem attained by reported methods in the literature.

The total emission obtained by the FIWCA is $0.2047$ ton per hour, while the total fuel cost achieved by the FIWCA is $\$941.8048$ per hour. Looking at Table \ref{Case5_comparison}, the FIWCA shows the best stability in obtained solutions, with the average emission of $0.2047$ ton per hour, compared with the other algorithms. The best obtained total emission values have been highlighted in bold in Table \ref{Case5_comparison}.

\subsection{IEEE 57-bus Test System}

\begin{figure}
 \centering
 \includegraphics[width=9cm, height=14cm,keepaspectratio]{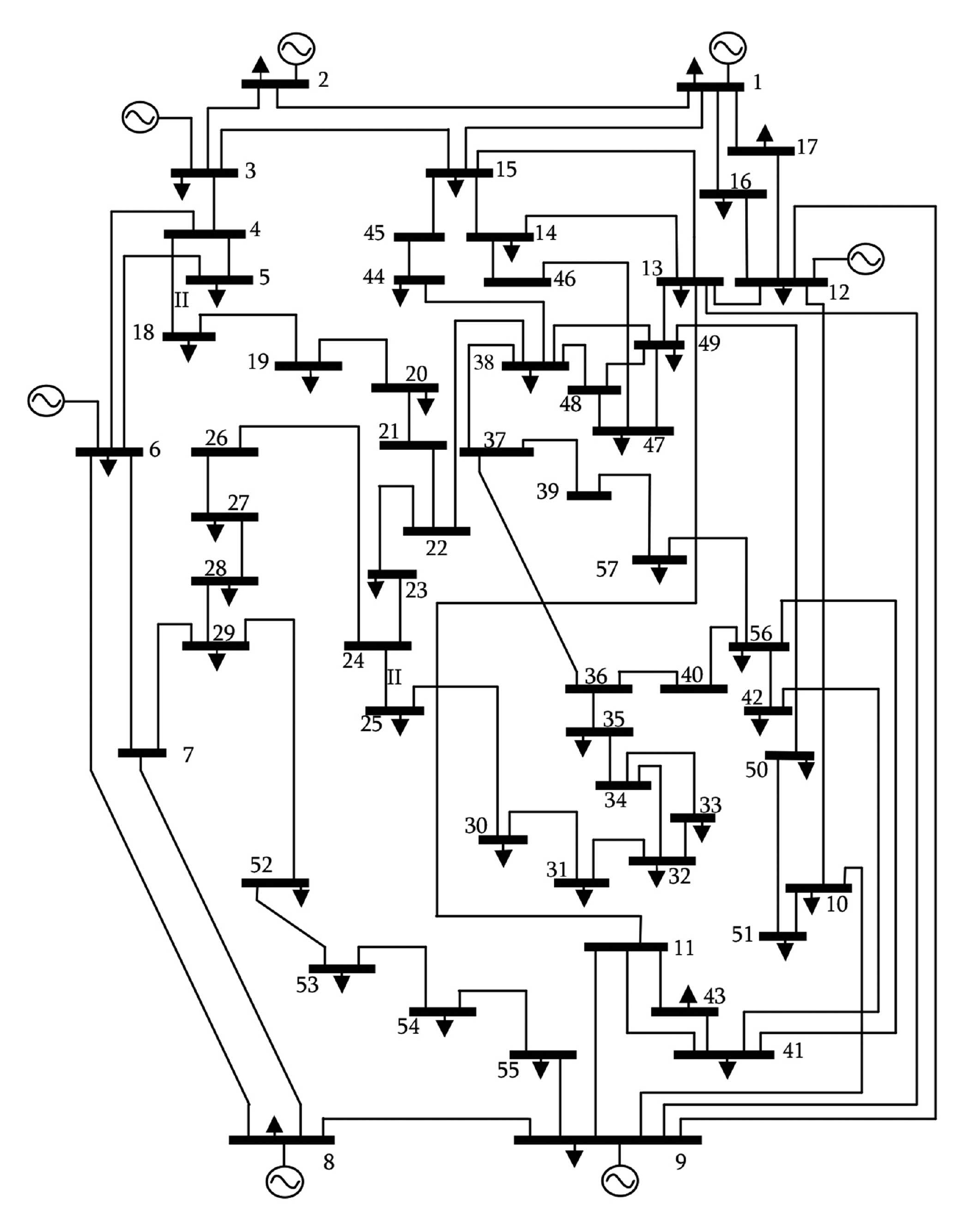}
 \caption{Schematic View of IEEE 57-bus Test System.}
 \label{Schematic_of_57_bus}
 \end{figure}

As mentioned in Section 1, one of the main features of the OPF problem, which makes it more complicated, is the size of the problem especially in large-scale power systems. To figure out the performance of the FIWCA in a larger system, the IEEE standard 57-bus test system has been considered. The IEEE standard 57-bus test system, as shown in Fig.~\ref{Schematic_of_57_bus}, consists of seven generators, 17 transformers with off-nominal tap ratio, and three shunt VAR compensation buses.
The lower and upper limits on control variables are available in Appendix, Table \ref{Case57_Upper_and_Lower_limits}. The WCA and FIWCA have been applied to this problem for a quadratic cost function given below
\begin{equation}
\label{FIWCA:Case57_cost}
J_{57-Bus} = \sum_{i=1}^{NG} f_i(P_{Gi}) = \sum_{i=1}^{NG} \big(c_{i0} + c_{i1} P_{Gi} + c_{i2} P_{Gi}^2 \big),
\end{equation}
where the coefficients of the \emph{i}-th generator $c_{i0}$, $c_{i1}$ and $c_{i2}$ are presented in Appendix, Table \ref{Case57:cost_coefficients}. The optimal control variables achieved by implementing the WCA and FIWCA for the IEEE 57-bus test system can be found in Table \ref{Case57_optimal_Values}.

\begin{table*}[!t]
\centering
\caption{The Best Control Design Variables Obtained by the FIWCA and WCA for IEEE 57-bus Test System.}\label{Case57_optimal_Values}
{\setlength\arrayrulewidth{1pt}
\begin{tabular}{llllll}
\hline
\multirow{2}{*}{\textbf{\begin{tabular}[c]{@{}l@{}}\textbf{Design} \\\textbf{Variables}\end{tabular}}} & \multirow{2}{*}{\textbf{WCA}} & \multirow{2}{*}{\textbf{FIWCA}} & \multirow{2}{*}{\textbf{\begin{tabular}[c]{@{}l@{}}\textbf{Design} \\\textbf{Variables}\end{tabular}}}& \multirow{2}{*}{\textbf{WCA}} & \multirow{2}{*}{\textbf{FIWCA}}\\  \\
\hline
$P_{G1}$  (MW)	& $143.7551$ & $142.9590$ & $T_{24-25}$ (p.u.)& $\hphantom{0}\hphantom{0}1.1000$ & $\hphantom{0}\hphantom{0}0.9977$  \\
$P_{G2}$  (MW)	& $\hphantom{0}90.3183$ & $\hphantom{0}91.4265$ & $T_{24-25}$  (p.u.)& $\hphantom{0}\hphantom{0}0.9691$ & $\hphantom{0}\hphantom{0}1.0486$  \\
$P_{G3}$  (MW)	& $\hphantom{0}44.8778$ & $\hphantom{0}45.0336$ &$ T_{24-26}$  (p.u.)& $\hphantom{0}\hphantom{0}1.0206$ & $\hphantom{0}\hphantom{0}1.0231$  \\
$P_{G6}$  (MW)	& $\hphantom{0}69.4730$ & $\hphantom{0}71.0449$ & $T_{7-29}$  (p.u.)& $\hphantom{0}\hphantom{0}0.9766$ & $\hphantom{0}\hphantom{0}0.9827$  \\
$P_{G8}$  (MW)	& $460.7114$ & $460.3489$ & $T_{34-32}$  (p.u.)& $\hphantom{0}\hphantom{0}0.9621$ & $\hphantom{0}\hphantom{0}0.9652$  \\
$P_{G9}$  (MW)	& $\hphantom{0}96.2347$ & $\hphantom{0}94.5440$ & $T_{11-41}$  (p.u.)& $\hphantom{0}\hphantom{0}0.9000$ & $\hphantom{0}\hphantom{0}0.9001$  \\
$P_{G12}$  (MW)	& $360.5641$ & $360.5006$ & $T_{15-45}$  (p.u.)& $\hphantom{0}\hphantom{0}0.9619$ & $\hphantom{0}\hphantom{0}0.9668$  \\
$V_{G1}$  (p.u.)	& $\hphantom{0}\hphantom{0}1.0600$ & $\hphantom{0}\hphantom{0}1.0597$ & $T_{14-46}$  (p.u.)& $\hphantom{0}\hphantom{0}0.9491$ & $\hphantom{0}\hphantom{0}0.9540$  \\
$V_{G2}$  (p.u.)	& $\hphantom{0}\hphantom{0}1.0558$ & $\hphantom{0}\hphantom{0}1.0574$ & $T_{10-51}$  (p.u.)& $\hphantom{0}\hphantom{0}0.9563$ & $\hphantom{0}\hphantom{0}0.9613$  \\
$V_{G3}$  (p.u.)	& $\hphantom{0}\hphantom{0}1.0440$ & $\hphantom{0}\hphantom{0}1.0494$ & $T_{13-49}$  (p.u.)& $\hphantom{0}\hphantom{0}0.9055$ & $\hphantom{0}\hphantom{0}0.9267$  \\
$V_{G6}$  (p.u.)	& $\hphantom{0}\hphantom{0}1.0545$ & $\hphantom{0}\hphantom{0}1.0566$ & $T_{11-43}$  (p.u.)& $\hphantom{0}\hphantom{0}0.9543$ & $\hphantom{0}\hphantom{0}0.9610$  \\
$V_{G8}$  (p.u.)	& $\hphantom{0}\hphantom{0}1.0589$ & $\hphantom{0}\hphantom{0}1.0600$ & $T_{40-56}$  (p.u.)& $\hphantom{0}\hphantom{0}0.9866$ & $\hphantom{0}\hphantom{0}0.9949$  \\
$V_{G9}$  (p.u.)	& $\hphantom{0}\hphantom{0}1.0326$ & $\hphantom{0}\hphantom{0}1.0362$ & $T_{39-57}$  (p.u.)& $\hphantom{0}\hphantom{0}0.9813$ & $\hphantom{0}\hphantom{0}0.9708$  \\
$V_{G12}$ (p.u.)	& $\hphantom{0}\hphantom{0}1.0314$ & $\hphantom{0}\hphantom{0}1.0393$ & $T_{9-55}$  (p.u.)& $\hphantom{0}\hphantom{0}0.9687$ & $\hphantom{0}\hphantom{0}0.9745$  \\
$T_{4-18}$  (p.u.)	& $\hphantom{0}\hphantom{0}1.0438$ & $\hphantom{0}\hphantom{0}0.9000$ & $Q_{C18}$  (MVAR)& $\hphantom{0}\hphantom{0}0.0699$ & $\hphantom{0}\hphantom{0}0.1289$  \\
$T_{4-18}$  (p.u.)	& $\hphantom{0}\hphantom{0}0.9313$ & $\hphantom{0}\hphantom{0}1.0792$ & $Q_{C25}$  (MVAR)& $\hphantom{0}\hphantom{0}0.1478$ & $\hphantom{0}\hphantom{0}0.1445$  \\
$T_{21-20}$  (p.u.)	& $\hphantom{0}\hphantom{0}1.0269$ & $\hphantom{0}\hphantom{0}1.0116$ & $Q_{C53}$  (MVAR)& $\hphantom{0}\hphantom{0}0.1318$ & $\hphantom{0}\hphantom{0}0.1312$  \\ \hline
\multirow{2}{*}{\textbf{\begin{tabular}[c]{@{}l@{}}\textbf{Fuel cost} \\\textbf{(\$/h)}\end{tabular}}}& \multirow{2}{*}{$41,678.9350$} & \multirow{2}{*}{$41,675.0794 $} & & \\ \\\hline
\multirow{2}{*}{\textbf{\begin{tabular}[c]{@{}l@{}}\textbf{Power loss} \\\textbf{((MW)}\end{tabular}}} & \multirow{2}{*}{$\hphantom{0}\hphantom{0}\hphantom{0}\hphantom{0}15.1347$} & \multirow{2}{*}{$\hphantom{0}\hphantom{0}\hphantom{0}\hphantom{0}15.0578$} & & \\ \\\hline
\end{tabular}}
\normalsize
\end{table*}

By evaluating the total optimized costs using the FIWCA, comparing with the previous works, as presented in Table \ref{Case57_comparison}, the superior performance of the FIWCA against the other considered optimizers is visible. The minimum value of the total generation fuel cost for this test system is $\$41,675.07$ per hour, best optimal solution obtained by different optimization algorithms. As demonstrated in Table \ref{Case57_optimal_Values}, not only the total fuel cost but also the total power loss obtained by the FIWCA is smaller than the standard WCA.

\begin{table}[]
\centering
\caption{Obtained Optimization Results for The IEEE 57-bus Test System, i.e. Fuel Costs (\$/h), Obtained by Several Optimizers.}\label{Case57_comparison}
{\setlength\arrayrulewidth{1pt}
\begin{tabular}{ll}
\hline
\textbf{Algorithms} & \multicolumn{1}{c}{\textbf{Fuel Cost}}  \\ \hline
\textbf{LDI-PSO \cite{adaryani2013artificial}}	 &	$41,815.50$ \\
\textbf{EADDE \cite{vaisakh2011evolving}}		 &	$41,713.62$ \\
\textbf{GSA \cite{duman2012optimal}	}	 &	$41,695.87$ \\
\textbf{TLBO \cite{ghasemi2015improved}} &	$41,695.66$ \\
\textbf{ABC \cite{adaryani2013artificial}}	 & $41,693.96$ \\ \hline
\textbf{FIWCA (This Study)}	 &	$\textbf{41,675.07}$ \\ \hline
\end{tabular}}
\normalsize
\end{table}

\subsection{Multi-period OPF problem with Renewable Resources}

Inserting the renewable energy resources into the power grid has many advantages, e.g., reducing the total network emission as the renewable sources are non-pollutant, requiring less maintainable comparing to the traditional power generation units, and decreasing the total power loss in transmission lines as the renewable resources are possible to be distributed within whole the network which decreases the amount of transmitted power through the lines. The total generated power by these renewable resources is dependent on the natural condition, i.e. wind speed, solar radiation, etc., which leads to time and location variant power output. On the other hand, the required power demand is variant during different periods of the day so it is required to solve a multi-period OPF problem.

The performance of FIWCA in solving the multi-period OPF problem including the renewable energy resources is the other challenge which would be investigated in this section. For this aim, the IEEE 30-Bus test system is considered for a period of 24 hours with 1-hour time steps with assumed variable load demands and distributed renewable resources with time and location variant power inserted at the load buses.

\begin{figure}[]
\centering
\subfloat[Demand Power and Generation Power]{\includegraphics[width=5.5cm, height=2.95cm]{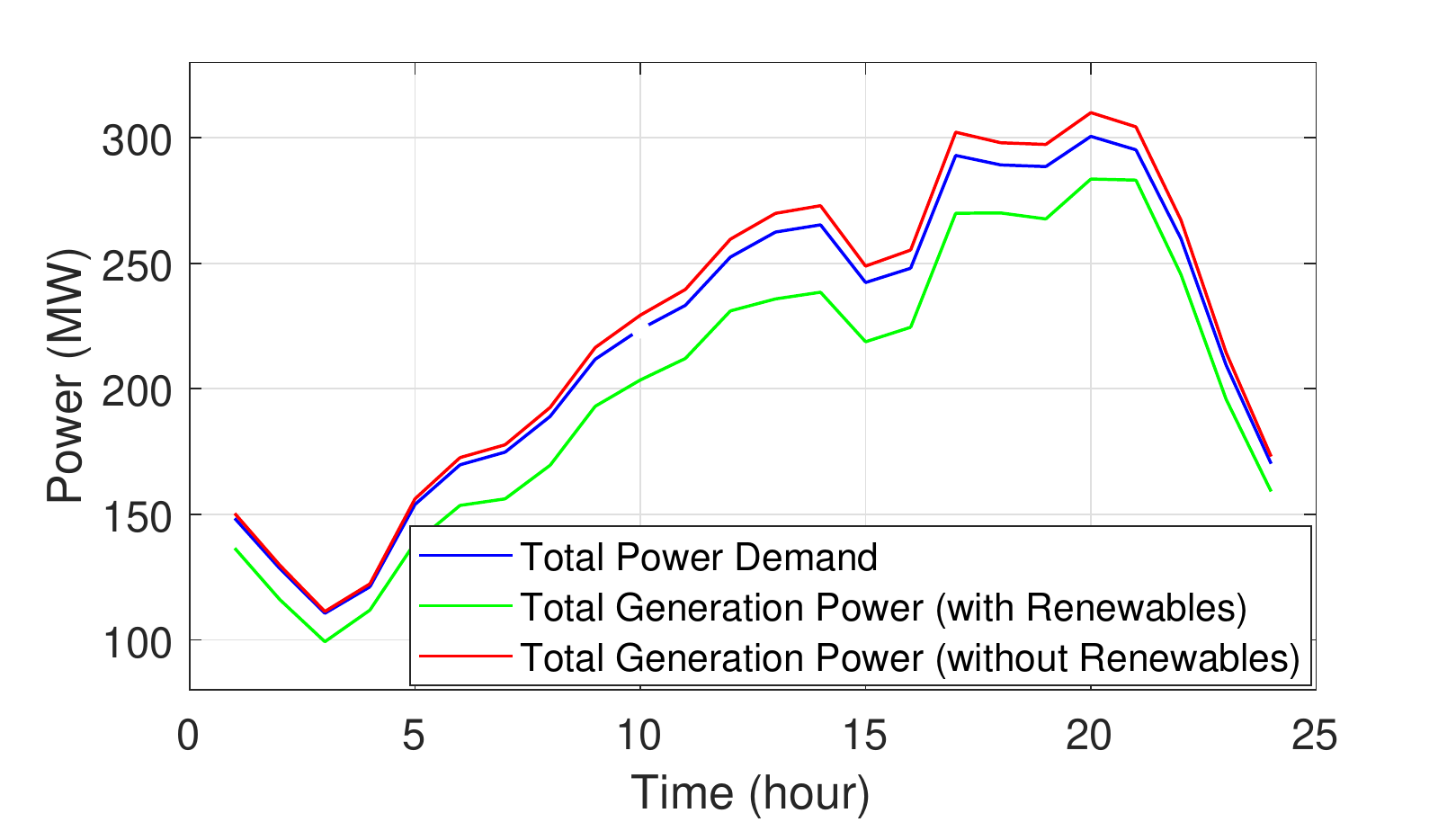}\label{ch3_fiwca24}} \hspace{0.5cm}
\subfloat[Difference Between Generated Power with and without Renewable Resources]{\includegraphics[width=5.5cm, height=8cm,keepaspectratio]{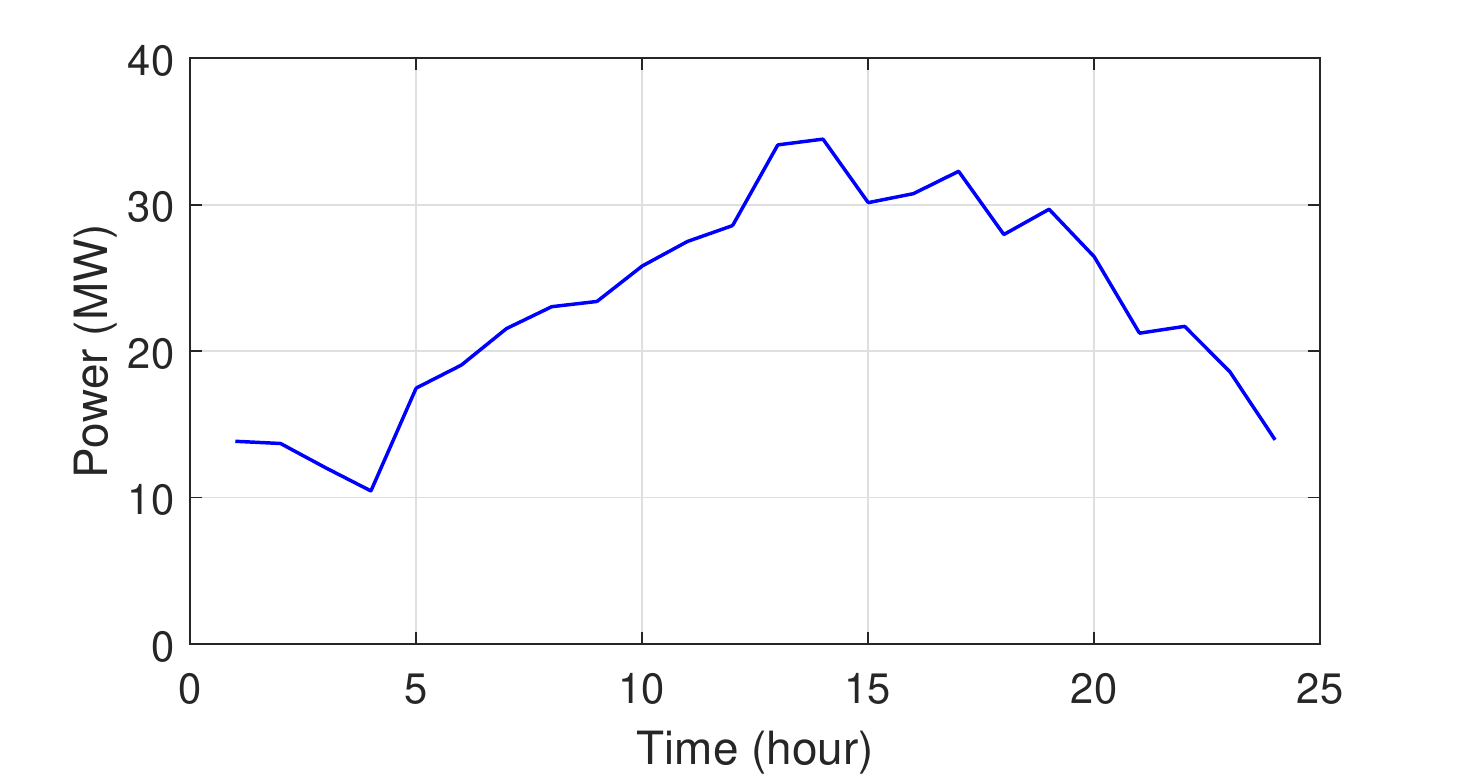}\label{ch3_gen_power}} \vspace{0.5cm}\\
\subfloat[Power Loss with and without Renewables]{\includegraphics[width=5.5cm, height=8cm,keepaspectratio]{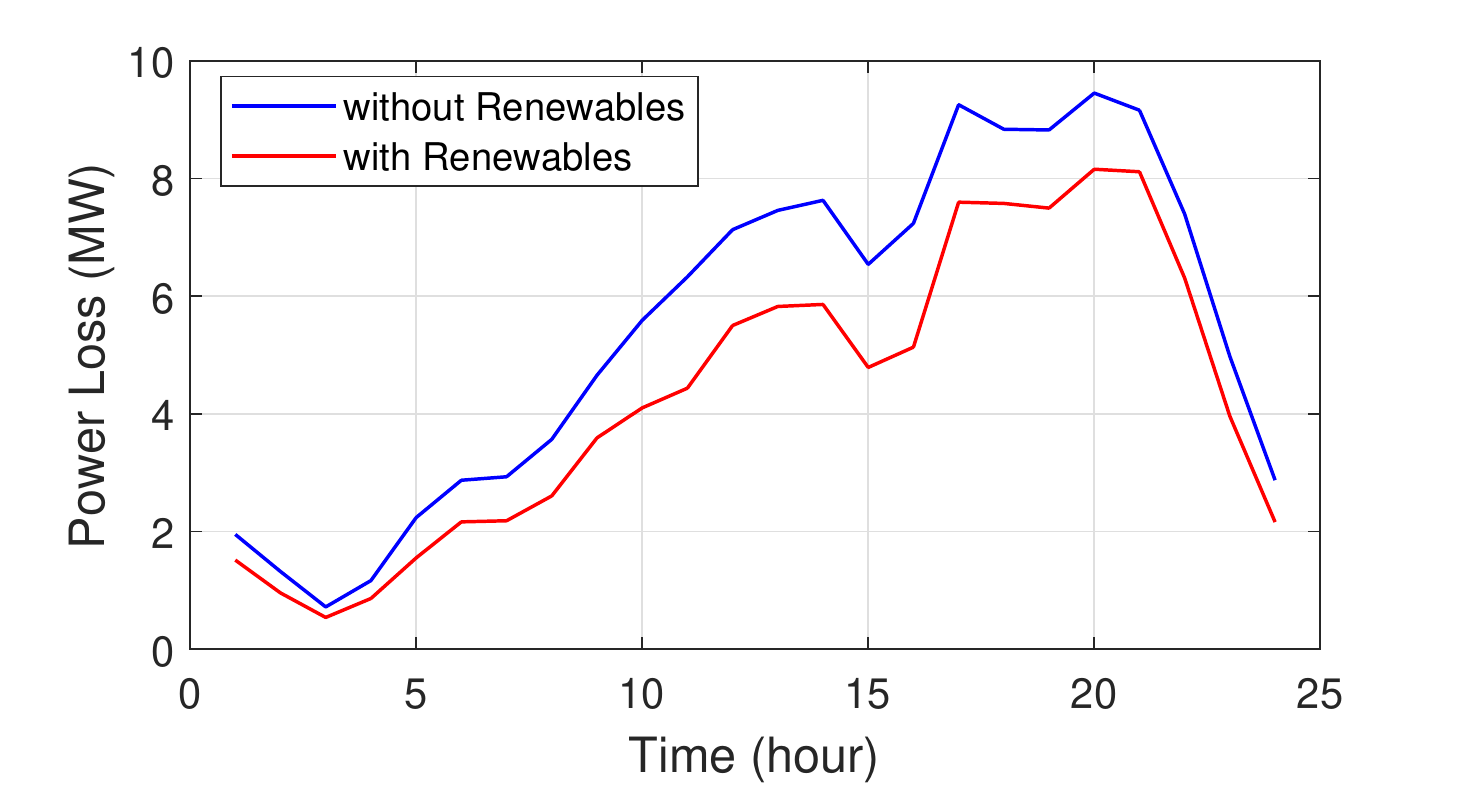}\label{ch3_power_loss}}\hspace{0.5cm}
\subfloat[Fuel Cost Savings by Implementing Renewable Resources]{\includegraphics[width=5.5cm, height=8cm,keepaspectratio]{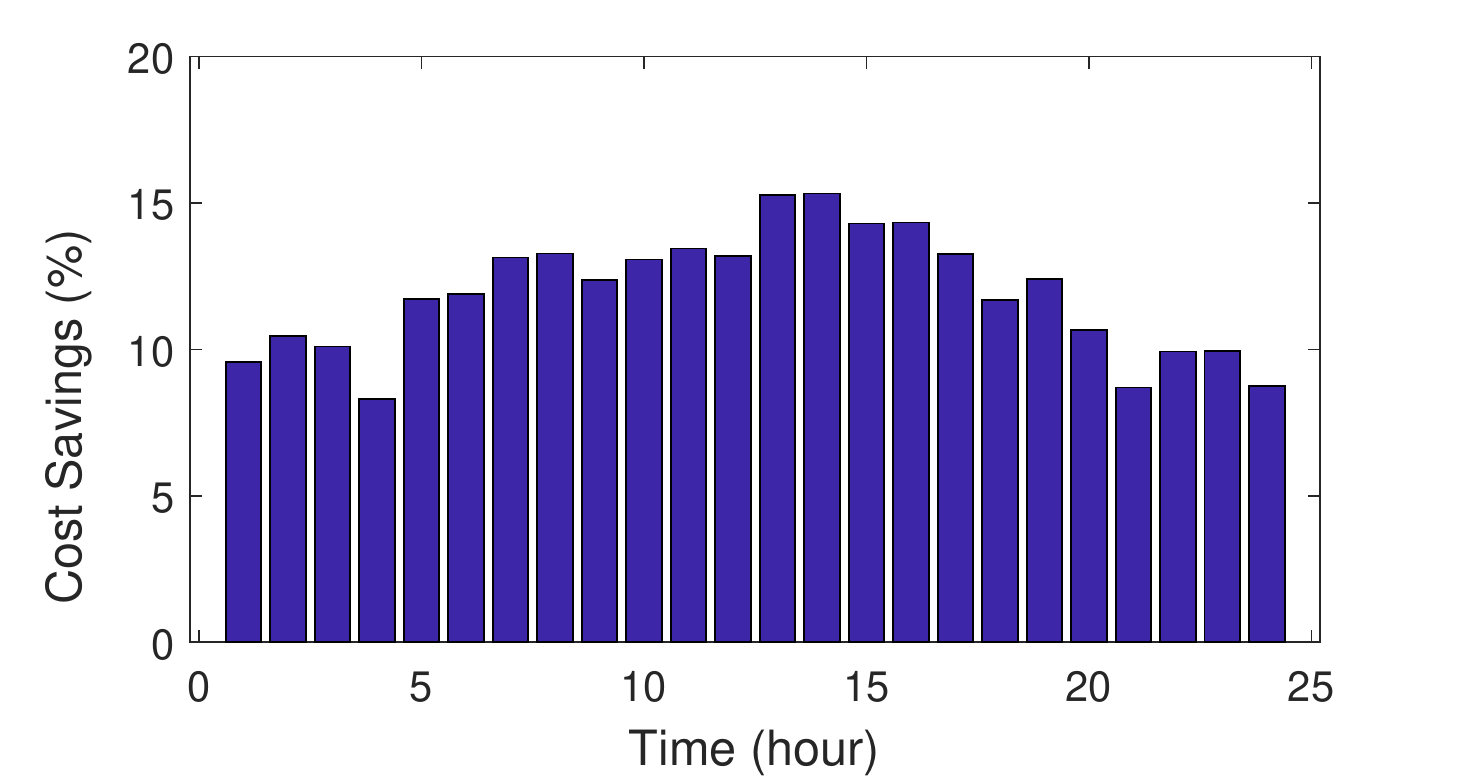}\label{ch3_cost_saving}}
\caption{Simulation Results of Solving Multi-period OPF Problem for IEEE 30-Bus Test System with and without Renewable Resources Using FIWCA}
\label{multi_period_OPF_renew}
\vspace*{1em}
\end{figure}
 To illustrate the performance of FIWCA in solving multi-period OPF problem with integrated renewable resources, Fig.~\ref{multi_period_OPF_renew} shows how to algorithm works in presence of renewable resources. Fig.~\ref{ch3_fiwca24} represents the power demand at each time step and total generated power by the traditional generators in two different scenarios, with integrated renewable resources and without them. The difference between the generated power of these scenarios is shown in Fig.~\ref{ch3_gen_power} which shows a total reduction of about $557.64$ $\$$ per day. The total power loss at each time is also demonstrated in Fig.~\ref{ch3_power_loss} and total cost saving percentage at each time period is presented in Fig~\ref{ch3_cost_saving} where the average saving during the day is about 10.4 $\%$.

 Table~\ref{multi_comparison} presents the comparison results of total fuel cost, total power loss and total emission for the period of 24 hours between the renewable integrated grid and the grid without any renewable resources. The simulation results show that total fuel cost is reduced about 12~$\%$, total power loss is approximately decreased 21~$\%$ and total emission from the fossil fuels is reduced around 7~$\%$.
\begin{table}[!t]
\centering
\caption{Comparison of Total Fuel Cost, Power Loss and Emission for Multi-period OPF problem with and without Renewable Resources Obtained by FIWCA.}
\label{multi_comparison}
\setlength\arrayrulewidth{1pt}
\begin{tabular}{llll}
\hline
\textbf{Test System} & \multicolumn{1}{l}{\textbf{\begin{tabular}[c]{@{}l@{}}Fuel Cost\\  (\$)\end{tabular}}} & \multicolumn{1}{l}{\textbf{\begin{tabular}[c]{@{}l@{}} Power Loss\\  (MW)  \end{tabular}}} & \multicolumn{1}{l}{\textbf{\begin{tabular}[c]{@{}l@{}}Emission \\ (tons)\end{tabular}}} \\ \hline
\textbf{IEEE 30-Bus without Renewables} &	$14111.73$ 	 & $ 130.10$ & $7.2178$\\
\textbf{IEEE 30-Bus with Renewables}	 &	$12402.43$   & 	$ 102.99$  & 	$ 6.6917$\\ \hline
\end{tabular}
\normalsize
\end{table}

\subsection{Discussion}
The simulation results obtained for the IEEE 30-bus and 57-bus test systems show the better performance of FIWCA in case of finding almost smaller cost function in different conditions, e.g. minimization of fuel cost, emission and voltage deviation compared with the standard WCA and the other reported optimizers. It is worth mentioning that the OPF problems are categorized as well-studied problems, whose solution quality is very tight to the globally optimal solution in some cases.

Therefore, the expected improvement level can be marginal, and no further major improvements can be performed. However, the statistical results validate the better quality of the proposed algorithm in term of stability of solution over replications which guarantees the reliability of the FIWCA comparing with the WCA and the other compared algorithms.

\begin{figure}[]
\vspace{0.5cm}
\centering
\centering
\subfloat[IEEE 30-bus Test Case System.]{\includegraphics[totalheight=0.17\textheight]{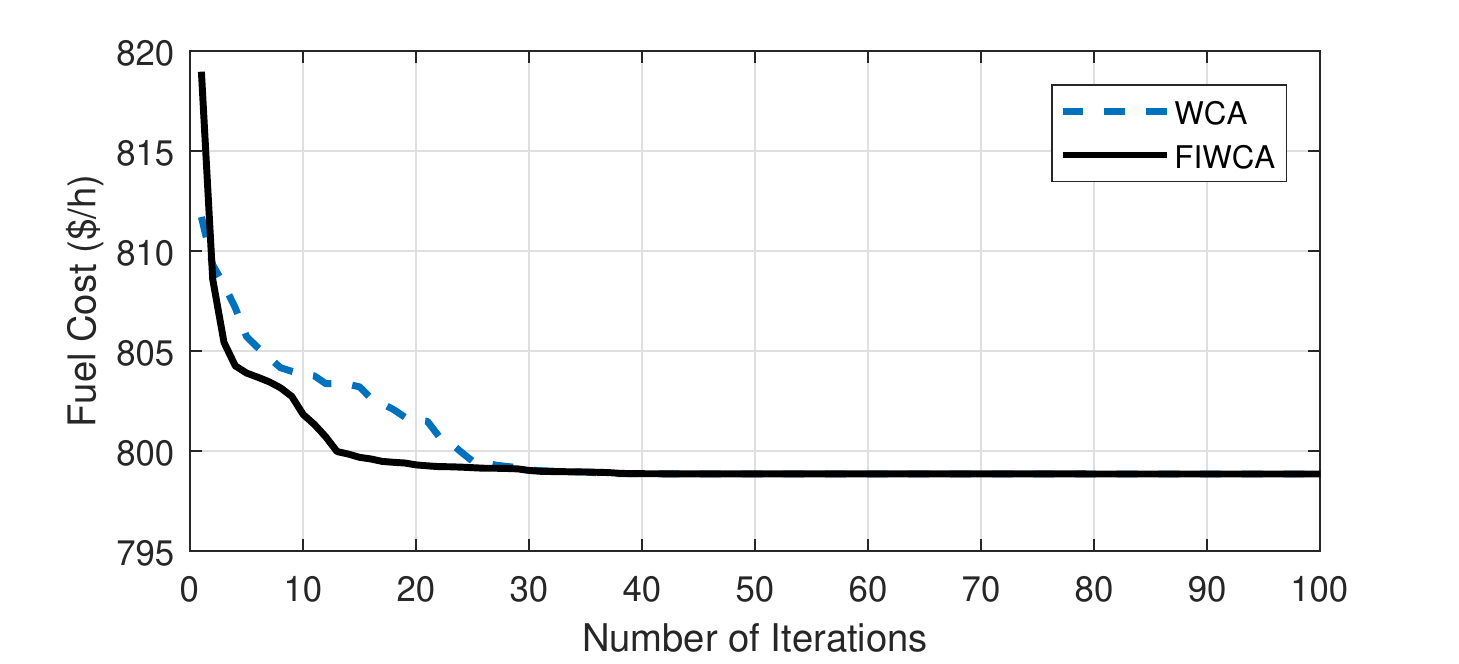}\label{30bus_conv}}\hspace{25mm}\\
\subfloat[IEEE 57-bus Test Case system.]{\includegraphics[totalheight=0.17\textheight]{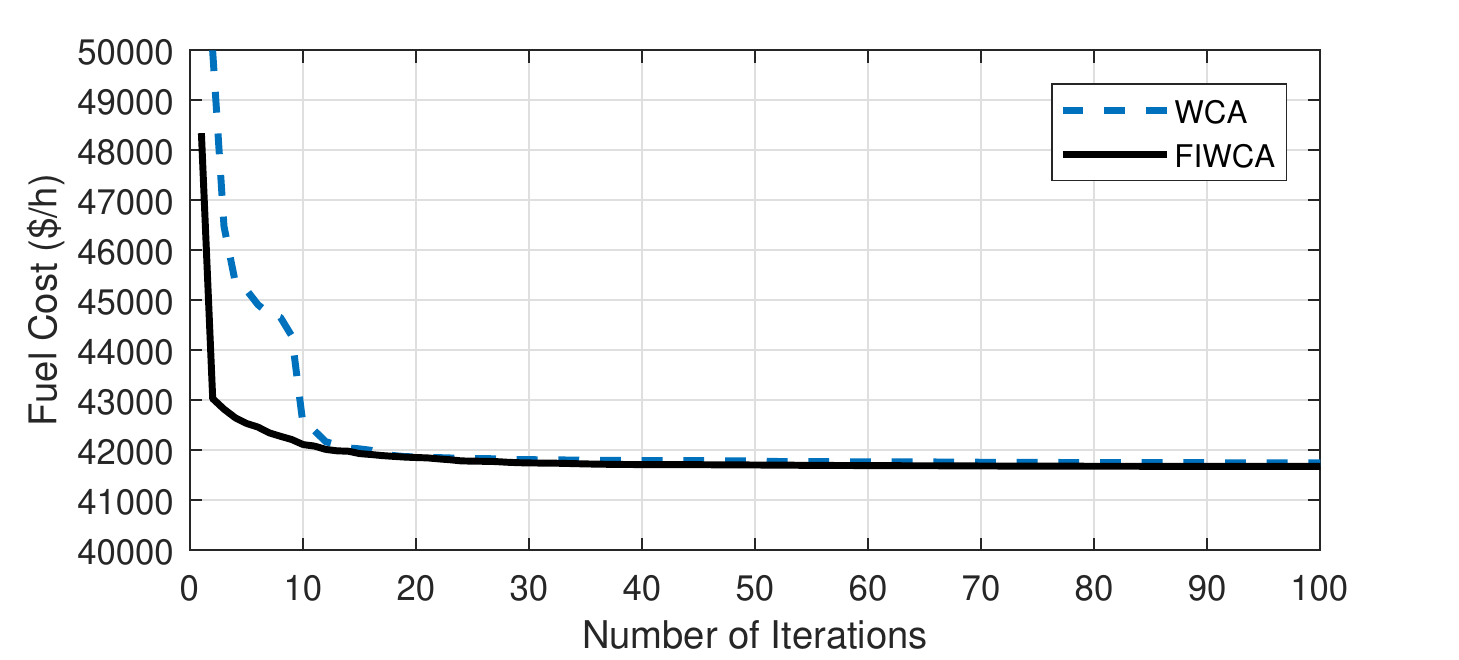}\label{57bus_conv}}\hspace{25mm}\\
\caption{Convergence Rate and Cost Reduction History Obtained by WCA and FIWCA.}
\label{Convergence_of_FIWCA}
\end{figure}

Furthermore, Fig.~\ref{Convergence_of_FIWCA} depicts the fuel cost value (\$/h) vs. iteration numbers for the IEEE 30-bus system, i.e. case 1, and IEEE 57-bus system which show the convergence characteristic of the proposed FIWCA. Looking at Fig.~\ref{Convergence_of_FIWCA}, the FIWCA has faster and mature convergence rate compared with the standard WCA, and it has converged within a reasonable number of iterations, i.e. around 20 iterations which is faster than the WCA with almost 30 iterations.

Table~\ref{ch3_comp_time} presents the comparison of computational time for solving the OPF problem obtained by different optimization algorithms. The optimal solutions and computational times show that the proposed FIWCA helps to enhance the quality of obtained solutions while being faster than the classical WCA and some other algorithms presented in Table~\ref{ch3_comp_time}.
\begin{table}[]
\vspace{5mm}
\centering
\caption{Comparison of Computational Time for The IEEE 30-bus and 57-bus Test Systems, Obtained by Different Optimizers.}
\label{ch3_comp_time}
\setlength\arrayrulewidth{1pt}
\begin{tabular}{lcc}
\hline
\multirow{2}{*}{\textbf{Algorithms}} & \multicolumn{2}{c}{\textbf{Computational Time (sec)}}       \\ \cline{2-3}
            & \multicolumn{1}{c}{\textbf{30-Bus}}     & \multicolumn{1}{c}{\textbf{57-Bus}}     \\ \hline
\textbf{DE-PS \cite{gitizadeh2014using}}	 &	$27.33$ &	N.A. \\
\textbf{TLBO \cite{ghasemi2015improved}}	 &	$24.27$ &	N.A. \\
\textbf{MDE \cite{sayah2008modified}}	     &	$23.25$ &	N.A.\\
\textbf{LTLBO \cite{ghasemi2015improved}}	 &	$22.14$ &	N.A. \\
\textbf{Fuzzy-GA \cite{hsiao2004optimal}}	 &	$22.07$ &	N.A. \\
\textbf{NPSO \cite{selvakumar2007new}}	     &	$20.45$ &	N.A. \\
\textbf{GSA \cite{adaryani2013artificial}}	 &	$19.71$ &	N.A.\\
\textbf{WCA (This Study) }	                 &	$20.38$ &	$48.70$ \\ \hline
\textbf{FIWCA (This Study)}	                 &	$\textbf{17.71}$ &	 $\textbf{43.35}$ \\ \hline
\end{tabular}
\normalsize
\end{table}

\section{Conclusions and future works}

This paper proposed an improved water cycle algorithm (WCA), called fully informed WCA (FIWCA) for solving the optimal power flow problem, a nonlinear and nonconvex optimization problem. The improvement to the WCA was the idea of using global and local information exchange for updating the positions of each individual in the population. Using fully informed algorithms help to increase the diversity in the selection of solution and consequently decrease the possibility of getting trapped at local optima, thereby enhancing the efficiency and accuracy significantly. The optimal power flow (OPF) problem studied different nonlinear constraints/objectives and evaluated the FIWCA for the optimality and efficiency of the algorithm. The OPF solution to IEEE 30-bus and 57-bus test systems was used to study the performance of the FIWCA with the existing approaches. Obtained optimization results demonstrated the computation benefits and efficiency of the FIWCA in reaching global optimal solutions with reasonable accuracy that is better than that of existing techniques such as standard WCA, genetic algorithm, particle swarm optimization, ant colony, and the others studied in the literature. Applying the proposed FIWCA to the multi-period OPF considering renewable energy fluctuations and storage devices is the future purpose of this investigation. Developing the proposed FIWCA to implement on the multi-objective OPF problem considering different objective values is another scope for future research works.



\bibliographystyle{IEEEtran}
\bibliography{FIWCA_2_Column}

\section{Appendix}\label{sec14}
The design variable limits and cost coefficients of different test cases for IEEE 30-bus and IEEE 57-bus systems are presented as following:
\begin{table}[ht]
\footnotesize
\centering
\caption{The Control Variables Limits for IEEE standard 30-bus \cite{lee1985united}.}\label{Case30_Upper_and_Lower_limits}
{\setlength\arrayrulewidth{1pt}
\begin{tabular}{llllll}
\hline
\textbf{Control Vars} &	\textbf{Min} &	\textbf{Max} &	\textbf{Control Vars} &	\textbf{Min} &		 \textbf{Max} \\ \hline
$P_{G1}$ (MW) & $50.00$ & $200.00$ & $T_{6-10}$  (p.u.) & $0.90$ & $1.10$  \\
$P_{G2}$  (MW) & $20.00$ & $\hphantom{0}80.00$ & $T_{4-12}$  (p.u.) & $0.90$ & $1.10$  \\
$P_{G5}$  (MW) & $15.00$ & $\hphantom{0}50.00$ & $T_{28-27}$  (p.u.) & $0.90$ & $1.10$  \\
$P_{G8}$  (MW) & $10.00$ & $\hphantom{0}35.00$ & $QC_{10}$  (MVAR) & $0.00$ & $5.00$  \\
$P_{G11}$  (MW) & $10.00$ & $\hphantom{0}30.00$ & $QC_{12}$  (MVAR) & $0.00$ & $5.00$  \\
$P_{G13}$  (MW) & $12.00$ & $\hphantom{0}40.00$ & $QC_{15}$  (MVAR) & $0.00$ & $5.00$  \\
$V_{G1}$  (p.u.) & $\hphantom{0}0.95$ & $\hphantom{0}\hphantom{0}1.10$ & $QC_{17}$  (MVAR) & $0.00$ & $5.00$  \\
$V_{G2}$  (p.u.) & $\hphantom{0}0.95$ & $\hphantom{0}\hphantom{0}1.10$ & $QC_{20}$  (MVAR) & $0.00$ & $5.00$  \\
$V_{G5}$  (p.u.) & $\hphantom{0}0.95$ & $\hphantom{0}\hphantom{0}1.10$ & $QC_{21}$  (MVAR) & $0.00$ & $5.00$  \\
$V_{G8}$  (p.u.) & $\hphantom{0}0.95$ & $\hphantom{0}\hphantom{0}1.10$ & $QC_{23}$  (MVAR)& $0.00$ & $5.00$  \\
$V_{G11}$  (p.u.) & $\hphantom{0}0.95$ & $\hphantom{0}\hphantom{0}1.10$ & $QC_{24}$  (MVAR) & $0.00$ & $5.00$  \\
$V_{G13}$  (p.u.) & $\hphantom{0}0.95$ & $\hphantom{0}\hphantom{0}1.10$ & $QC_{29}$  (MVAR) & $0.00$ & $5.00$  \\
$T_{6-9}$  (p.u.) & $\hphantom{0}0.90$ & $\hphantom{0}\hphantom{0}1.10$	& &  &\\ \hline		
\end{tabular}}
\normalsize
\end{table}
\begin{table}[ht]
\footnotesize
\centering
\caption{Generation Cost Coefficients for Case 1 (IEEE 30-bus) \cite{lee1985united}.}\label{Case1:cost_coefficients}
{\setlength\arrayrulewidth{1pt}
\begin{tabular}{llll}
\hline
\multirow{2}{*}{\textbf{Bus No.}} & \multicolumn{3}{c}{\textbf{Cost Coefficients}}       \\ 
            & \multicolumn{1}{c}{$\mathbf{c_{i0}}$}     & \multicolumn{1}{c}{$\mathbf{c_{i1}}$}    & \multicolumn{1}{c}{$\mathbf{c_{i2}}$}  \\ \hline
\textbf{\hphantom{0}1}	& $0.00$	& $2.00$	& $0.00375$ \\
\textbf{\hphantom{0}2}	& $0.00$	& $1.75$	& $0.01750$ \\
\textbf{\hphantom{0}5}	& $0.00$	& $1.00$	& $0.06250$ \\
\textbf{\hphantom{0}8}	& $0.00$	& $3.25$	& $0.00834$ \\
\textbf{11}	& $0.00$	& $3.00$	& $0.02500$ \\
\textbf{13}	& $0.00$	& $3.00$	& $0.02500$ \\ \hline
\end{tabular}}
\normalsize
\end{table}
\begin{table}[ht]
\footnotesize
\centering
\caption{Generation Cost Coefficients for Case 3 (IEEE 30-bus) \cite{abido2002optimal}.}\label{Case3:cost_coefficients}
{\setlength\arrayrulewidth{1pt}
\begin{tabular}{llllll}
\hline
\multirow{2}{*}{\textbf{Bus No.}} & \multirow{2}{*}{\textbf{From (MW)}} & \multirow{2}{*}{\textbf{To (MW)}} & \multicolumn{3}{c}{\textbf{Cost Coefficients}}                        \\ 
     & &       & \multicolumn{1}{c}{$\mathbf{c_{i0_k}}$}     & \multicolumn{1}{c}{$\mathbf{c_{i1_k}}$}    & \multicolumn{1}{c}{$\mathbf{c_{i2_k}}$}   \\ \hline
\multirow{2}{*}{\textbf{1}} & $\hphantom{0}50$  & $140$  & $55.00$	& $0.70$	& $0.0050$ \\
                            & $140$  & $200$  	& $82.50$	& $1.05$	& $0.0075$ \\
\multirow{2}{*}{\textbf{2}} & $\hphantom{0}20$  & $\hphantom{0}55$ 	& $40.00$	& $0.30$	& $0.0100$ \\
                                 &  $\hphantom{0}55$  & $\hphantom{0}80$  	& $80.00$	& $0.60$	& $0.2000$ \\\hline
\end{tabular}}
\normalsize
\end{table}
\begin{table}[ht]
\footnotesize
\centering
\caption{Generation Cost Coefficients for Case 4 (IEEE 30-bus) \cite{selvakumar2007new}.}\label{Case4:cost_coefficients}
{\setlength\arrayrulewidth{1pt}
\begin{tabular}{llllllll}
\hline
\multirow{2}{*}{\textbf{Bus No.}} & \multirow{2}{*}{\textbf{\begin{tabular}[c]{@{}l@{}}$\mathbf{P_{Gi}^{min}}$ \\ (MW)\end{tabular}}}& \multirow{2}{*}{\textbf{\begin{tabular}[c]{@{}l@{}}$\mathbf{P_{Gi}^{max}}$ \\ (MW)\end{tabular}}} & \multicolumn{5}{c}{\textbf{Cost Coefficients}}                        \\ 
            &   &       & \multicolumn{1}{c}{$\mathbf{c_{k0}}$}     & \multicolumn{1}{c}{$\mathbf{c_{k1}}$}    & \multicolumn{1}{c}{$\mathbf{c_{k2}}$}   & \multicolumn{1}{c}{$\mathbf{d_k}$} & \multicolumn{1}{c}{$\mathbf{e_k}$}\\ \hline
\textbf{1} & \multicolumn{1}{c}{$50$}  & \multicolumn{1}{c}{$200$}  & $150$  & $2.00$  & $0.0016$  & $50$  & $0.0630$ \\
\textbf{2} & \multicolumn{1}{c}{$20$}  & \multicolumn{1}{c}{$\hphantom{0}80$}  & $\hphantom{0}25$  & $2.50$  & $0.0100$  & $40$  & $0.0980$ \\ \hline
\end{tabular}}
\normalsize
\end{table}
\begin{table}[ht]
\footnotesize
\centering
\caption{Generation Emission Coefficients used in Case 5 (IEEE 30-bus) \cite{lee1985united}.}\label{Case5:cost_coefficients}
{\setlength\arrayrulewidth{1pt}
\begin{tabular}{llllll}
\hline
\multirow{2}{*}{\textbf{Bus No.}}   & \multicolumn{5}{c}{\textbf{Cost Coefficients}}       \\ 
         & \multicolumn{1}{c}{$\mathbf{c_{k0}^{SOX}}$}     & \multicolumn{1}{c}{$\mathbf{c_{k1}^{SOX}}$}   & \multicolumn{1}{c}{$\mathbf{c_{k2}^{SOX}}$}  & \multicolumn{1}{c}{$\mathbf{d_{k}^{NOX}}$} & \multicolumn{1}{c}{$\mathbf{e_{k}^{NOX}}$}\\ \hline
\textbf{\hphantom{0}1} &	$0.04091$  & $-0.05554$  & $0.06490$  & $0.000200$  & $2.857$ \\
\textbf{\hphantom{0}2}	 & $0.02543$  & $-0.06047$  & $0.05638$  & $0.000500$  & $3.333$  \\
\textbf{\hphantom{0}5}	& $0.04258$  & $-0.05094$  & $0.04586$  & $0.000001$  & $8.000$ \\
\textbf{\hphantom{0}8}	& $0.05326$  & $-0.03550$  & $0.03380$  & $0.002000$  & $2.000$ \\
\textbf{11} &	$0.04258$  & $-0.05094$  & $0.04586$  & $0.000001$  & $8.000$ \\
\textbf{13}	& $0.06131$  & $-0.05555$  & $0.05151$  & $0.000010$  & $6.667$ \\ \hline
\end{tabular}}
\normalsize
\end{table}
\begin{table}[ht]
\footnotesize
\centering
\caption{The Control Variables Limits for the IEEE 57-bus Test System \cite{vaisakh2011evolving}.}\label{Case57_Upper_and_Lower_limits}
{\setlength\arrayrulewidth{1pt}
\begin{tabular}{llllll}
\hline
\textbf{Control Vars} &	\textbf{Min} &	\textbf{Max} &	\textbf{Control Varss} &	\textbf{Min} &		 \textbf{Max} \\ \hline
$P_{G1}$ (MW) & $0.00$ & $575.88$ & $T_{24-25}$ (p.u.) & $0.90$ & $\hphantom{0}1.10$  \\
$P_{G3}$ (MW) & $0.00$ & $140.00$ & $T_{24-26}$ (p.u.) & $0.90$ & $\hphantom{0}1.10$  \\
$P_{G6}$ (MW) & $0.00$ & $100.00$ & $T_{7-29}$ (p.u.) & $0.90$ & $\hphantom{0}1.10$  \\
$P_{G8}$ (MW) & $0.00$ & $550.00$ & $T_{34-32}$ (p.u.) & $0.90$ & $\hphantom{0}1.10$  \\
$P_{G9}$ (MW) & $0.00$ & $100.00$ & $T_{11-41}$ (p.u.) & $0.90$ & $\hphantom{0}1.10$  \\
$P_{G12}$ (MW) & $0.00$ & $410.00$ & {$T_{15-45}$ (p.u.)} & $0.90$ & $\hphantom{0}1.10$  \\
$V_{G1}$ (p.u.) & $0.94$ & $\hphantom{0}\hphantom{0}1.06$ & $T_{14-46}$ (p.u.) & $0.90$ & $\hphantom{0}1.10$  \\
$V_{G2}$ (p.u.) & $0.94$ & $\hphantom{0}\hphantom{0}1.06$ & $T_{10-51}$ (p.u.) & $0.90$ & $\hphantom{0}1.10$  \\
$V_{G3}$ (p.u.) & $0.94$ & $\hphantom{0}\hphantom{0}1.06$ & $T_{13-49}$ (p.u.) & $0.90$ & $\hphantom{0}1.10$  \\
$V_{G6}$ (p.u.) & $0.94$ & $\hphantom{0}\hphantom{0}1.06$ & $T_{11-43}$ (p.u.) & $0.90$ & $\hphantom{0}1.10$  \\
$V_{G8}$ (p.u.) & $0.94$ & $\hphantom{0}\hphantom{0}1.06$ & $T_{40-56}$ (p.u.) & $0.90$ & $\hphantom{0}1.10$  \\
$V_{G9}$ (p.u.) & $0.94$ & $\hphantom{0}\hphantom{0}1.06$ & $T_{39-57}$ (p.u.)  & $0.90$ & $\hphantom{0}1.10$  \\
$V_{G12}$ (p.u.) & $0.94$ & $\hphantom{0}\hphantom{0}1.06$ & $T_{9-55}$ (p.u.) & $0.90$ & $\hphantom{0}1.10$  \\
$T_{4-18}$ (p.u.) & $0.90$ & $\hphantom{0}\hphantom{0}1.10$ &$QC_{18}$ (MVAR) & $0.00$ & $30.00$  \\
$T_{4-18}$ (p.u.) & $0.90$ & $\hphantom{0}\hphantom{0}1.10$ & $QC_{25}$ (MVAR) & $0.00$ & $30.00$  \\
$T_{21-20}$ (p.u.) & $0.90$ & $\hphantom{0}\hphantom{0}1.10$ & $QC_{53}$ (MVAR) & $0.00$ & $30.00$ \\ \hline
\end{tabular}}
\normalsize
\end{table}
\begin{table}[ht]
\footnotesize
\centering
\caption{Generation Cost Coefficients for IEEE 57-bus \cite{vaisakh2011evolving}.}\label{Case57:cost_coefficients}
{\setlength\arrayrulewidth{1pt}
\begin{tabular}{llll}
\hline
\multirow{2}{*}{\textbf{Bus No.}} & \multicolumn{3}{c}{\textbf{Cost Coefficients}}       \\ 
            & \multicolumn{1}{c}{$\mathbf{c_{k0}}$}     & \multicolumn{1}{c}{$\mathbf{c_{k1}}$}    & \multicolumn{1}{c}{$\mathbf{c_{k2}}$}   \\ \hline
\textbf{1}	& $0.00$  & $20.00$  & $0.0775$ \\
\textbf{2}	& $0.00$  & $40.00$  & $0.0100$ \\
\textbf{3}	& $0.00$  & $20.00$  & $0.2500$ \\
\textbf{6}	& $0.00$  & $40.00$  & $0.0100$ \\
\textbf{8}	& $0.00$  & $20.00$  & $0.0222$ \\
\textbf{9}	& $0.00$  & $40.00$  & $0.0100$ \\ \hline
\end{tabular}}
\normalsize
\end{table}

\end{document}